\documentclass[sn-mathphys,Numbered]{sn-jnl}% Math and Physical Sciences Reference Style
\usepackage{graphicx}%
\usepackage{multirow}%
\usepackage{amsmath,amssymb,amsfonts}%
\usepackage{amsthm}%
\usepackage{mathrsfs}%
\usepackage[title]{appendix}%
\usepackage{xcolor}%
\usepackage{textcomp}%
\usepackage{manyfoot}%
\usepackage{booktabs}%
\usepackage{algorithm}%
\usepackage{algorithmicx}%
\usepackage{algpseudocode}%
\usepackage{listings}%
\usepackage{graphicx}
\usepackage{color}
\usepackage{pifont}
\usepackage{subfigure}
\usepackage[section]{placeins}
\usepackage{makecell}
\usepackage{appendix}
\theoremstyle{thmstyleone}%
\newtheorem{theorem}{Theorem}%  meant for continuous numbers
\newtheorem{lemma}{Lemma}
\theoremstyle{thmstyletwo}%

\theoremstyle{thmstylethree}%
\newtheorem{definition}{Definition}%
\newtheorem{corollary}{Corollary}

\numberwithin{equation}{section}
\numberwithin{figure}{section}
\numberwithin{table}{section}

\begin{document}

\title[Stochastic Conjugate Frameworks for Nonconvex and Nonsmooth Optimization]{Stochastic Conjugate Frameworks for Nonconvex and Nonsmooth Optimization}

\author[]{\fnm{Jiangshan} \sur{Wang}}\email{wjs\_13977326105@163.com}
\author[]{\fnm{Zheng} \sur{Peng}}\email{pzheng@xtu.edu.cn}


\affil[]{\orgdiv{School of Mathematics and Computational Science}, \orgname{Xiangtan University}, \orgaddress{\city{Xiangtan}, \postcode{411105}, \state{Hunan Province}, \country{P. R. China}}}

%\affil[2]{\orgdiv{Department}, \orgname{Organization}, \orgaddress{\street{Street}, \city{City}, \postcode{10587}, \state{State}, \country{Country}}}

%\affil[3]{\orgdiv{Department}, \orgname{Organization}, \orgaddress{\street{Street}, \city{City}, \postcode{610101}, \state{State}, \country{Country}}}

\abstract{We introduce two new stochastic conjugate frameworks for a class of nonconvex and possibly also nonsmooth optimization problems. These frameworks are built upon Stochastic Recursive Gradient Algorithm (SARAH) and we thus refer to them as Acc-Prox-CG-SARAH and Acc-Prox-CG-SARAH-RS, respectively. They are efficiently accelerated, easy to implement, tune free and can be smoothly extended and modified. We devise a deterministic restart scheme for stochastic optimization and apply it in our second stochastic conjugate framework, which serves the key difference between the two approaches. In addition, we apply the ProbAbilistic Gradient Estimator (PAGE) and further develop a practical variant, denoted as Acc-Prox-CG-SARAH-ST, in order to reduce potential computational overhead. We provide comprehensive and rigorous convergence analysis for all three approaches and establish linear convergence rates for unconstrained minimization problem with nonconvex and nonsmooth objective functions. Experiments have demonstrated that Acc-Prox-CG-SARAH and Acc-Prox-CG-SARAH-RS both outperform state-of-art methods consistently and Acc-Prox-CG-SARAH-ST can as well achieve comparable convergence speed. In terms of theory and experiments, we verify the strong computational efficiency of the deterministic restart scheme in stochastic optimization methods.}


\maketitle

\section{Introduction}
In the paper, we are concerned with minimizing the following composite, nonconvex, and possibly also nonsmooth optimization problem:
\begin{equation}
\label{problem}
\min _{w \in \mathbb{R}^d}\left\{P(w)=f(w)+\varphi(w) = \frac{1}{n} \sum_{i=1}^n f_i(w)+\varphi(w)\right\},
\end{equation}
where $f$ is the average of $n$ smooth functions $f_1,\ldots,f_n$ from $\mathbb{R}^d$ to $\mathbb{R}$, and $\varphi$: $\mathbb{R}^d\to\mathbb{R}\cup\{+\infty\}$ is a proper,
closed, and convex function, possibly non-smooth. If $\varphi$ is a relatively simple and convex regularizer, problem (\ref{problem}) arises in a wide range of applications in machine learning and statistics, such as image processing (Sitzmann et al., 2018 \cite{sitzmann2018image}), neural networks (Bottou et al., 2018 \cite{bottou2018neural}) and parallel deep learning (Guo et al., 2020 \cite{guo2020deep}). Further, if $\varphi$ is the indicator of a nonempty, closed, and convex set, then (\ref{problem}) also covers constrained nonconvex optimization problems.   

\subsection{Stochastic First-Order Methods}
Stochastic first-order methods have surged into prominence in large-scale optimization, which are of three leading categories with the stochastic gradient-types, conditional
gradient methods and primal-dual approaches. 

The stochastic gradient-types mainly rely on the proximal stochastic gradient descent (PSGD) (Ghadimi and Lan, 2013 \cite{ghadimi2013proximal}), which can be described by
$$
w_{t+1}=\arg\min_{w\in\mathbb{R}^{d}}\left\{\frac{1}{2\eta_{t}}\|w-w_{t}\|^{2}+\langle w,  g(w_{t}, \nu_t)\rangle+\varphi(w)\right\},
$$
where $g$ is a stochastic estimate of $\nabla f$, and $\nu_t$ is a random variable. It becomes projection gradient descent when $\varphi$ is the indicator function. Prevalent choices for $g$ involve vanilla stochastic gradient (SGD) (Robbins and Monro, 1951 \cite{robbins1951SGD}), stochastic average gradient (SAG/SAGA) (Roux et al., 2012 \cite{roux2012SAG}, Defazio et al., 2014 \cite{defazio2014SAGA}),  stochastic variance reduced gradient (SVRG) (Johnson and Zhang, 2013 \cite{johnson2013SVRG}), semi-stochastic gradient (S2GD) (Kone{\v{c}}n{\`y} et al. 2017 \cite{konevcny2017S2GD}) and stochastic recursive gradient (SARAH) (Nguyen et al. 2017 \cite{nguyen2017SARAH}), etc. Related proximal variants encompass ProxSAG (Schmidt et al., 2017 \cite{schmidt2017proxSAG}), ProxSVRG+ (Li and Li, 2018 \cite{li2018proxsvrg+}), PROXASAGA (Pedregosa et al., 2017 \cite{pedregosa2017proxASAGA}) and proxSARAH (Pham et al., 2020 \cite{phamproxsarah}), etc. In contrast, the computational cost of PSGD corresponds to $1/n$ that of the PGD per iteration.

In terms of variance reduction, mini-batch and importance sampling schemes are intensely studied in both convex and nonconvex settings (see, e.g., Ghadimi et al. (2016) \cite{ghadimi2016batch} for nonconvex optimization and Liu et al. (2020) \cite{liu2020SARAHI} under convexity setting). Besides, several algorithms have shown potential to speed up variance reduction. The remarkable variants include SVRG-based approaches, e.g., Katyusha (Allen-Zhu, 2017 \cite{allen2017katyusha}), MiG (Zhou et al., 2018 \cite{zhou2018MiG}), SNVRG (Zhou et al., 2018 \cite{zhou2018SNVRG}); SARAH-related variants, e.g., SPIDERE (Fang et al., 2018 \cite{fang2018spider}), SpiderBoost (Wang et al., 2019 \cite{wang2019spiderboost}); and others, e.g., SEGA (Hanzely et al., 2018 \cite{hanzely2018SEGA}), Prox-SDCA (Zhao and Zhang, 2015 \cite{zhao2015proxSDCA}), etc. In general, they often require stronger assumptions.

Techniques have been proposed to speed up convergence, and the following are the most notable. The first class relies in Nesterov's acceleration, see, e.g., Murata and Suzuki (2017) \cite{murata2017nesterov} and Lan and Zhou (2018) \cite{lan2018nesterov}. The second one is based on the choice of growing epoch length, including Allen-Zhu and Yang (2016) \cite{allen2016epochlength}. The third idea utilizes momentum acceleration, such as ARMD (Hien et al., 2018 \cite{hien2019momentum}), which combines accelerated proximal gradient descent (AGD) with the inexact computation of proximal points. Also, Li et al. (2020) \cite{li2020momentum} provided a comprehensive survey of various momentum schemes to achieve accelerated algorithms. Typically, accelerations require more complicated algorithmic designs and parameter tunings.

Stochastic first-order methods also include primal-dual approach, e.g., the stochastic primal-dual hybrid gradient algorithm (SPDHG) (Chambolle et al., 2018 \cite{chambolle2018primaldual}). The conditional gradient methods involve stochastic Frank-Wolfe (SFW) approach (Reddi et al., 2016 \cite{reddi2016SFW}), which is efficient in nonconvex stochastic optimization problems.

\subsection{Stochastic Conjugate Gradient (SCG) Methods}
Conjugate gradient (CG) method comprises a class of unconstrained optimization algorithms with low storage requirements and strong convergence properties. It is generally believed to be one of the most efficient methods for solving unconstrained optimization problems, particularly those on a large scale.

In terms of applications, the SCG methods have been intensely explored in stochastic optimization for machine learning problems. Liu et al. (2018) \cite{liu2018CG1} developed a class of stochastic conjugate gradient descent (SCGDM) algorithms to optimize the weight parameters in the neural language model. Huang and Zhou (2015) \cite{huang2015CG2} applied the stochastic conjugate gradient (SCG) approach to least-squares seismic inversion problems for function approximation. And Li et al. (2018) \cite{li2018CG3} derived a preconditioned stochastic conjugate gradient (PSCG) method for the image processing problems. Wei et al. (2020) \cite{wei2020CG5} introduced a learned conjugate gradient descent network (LcgNet) and utilized it for massive multiple-input multiple-output (MIMO) detection.

Due to the noisy gradient estimates, vanilla SCG method has suffered low convergence rate or divergence results. Therefore, several improved SCG-type algorithms have been presented. Based on variance reduced estimator SVRG, Jin et al. (2018) \cite{jin2018cg-svrg} proposed CGVR algorithm, which exhibits faster convergence than SVRG. Afterward, Xue et al. (2021) \cite{xue2021CG4} proposed an online SIFR-CG algorithm, which combined the SVRG estimator and improved Fletcher-Reeves (IFR) approach. It achieves a linear convergence rate under strong Wolfe line search in strongly convex cases. To reduce computational workload, Kou and Yang (2022) \cite{kou2022cgsaga} developed SCGA algorithm upon the SAGA estimator. SCGA obtains a linear rate of convergence for strongly convex functions and consumes less computation cost than CGVR. The SCG method has also been developed based on biased estimators, such as SARAH. Yang (2022) \cite{yang2022cgsarahso} \cite{yang2022cgsarahhd} proposed CG-SARAH algorithm and combine it with the random Barzilai and Borwein (RBB) step size (Yang et al, 2018 \cite{yang2018random}) and hyper-gradient descent (HD) approach (Baydin et al., 2017 \cite{baydin2017HD}), respectively. Several practical variants have also been developed to work with an automatic and adaptive selection for the epoch length. Inspired by ordinary differential equations (ODEs), Yang (2023) \cite{yang2023PBCG} developed adaptive powerball stochastic conjugate gradient algorithm (pbSGD) for nonconvex optimization, the HD approach was introduced into pbSGD to obtain an online learning rate.

\subsection{Motivation}
Due to significant applications in deep learning, large-scale nonconvex problems have garnered immense attention in recent years. Hitherto, most state-of-the-art numerical algorithms focus on solving the nonconvex problems in non-composite setting (see, e.g., Fang et al., 2018 \cite{fang2018spider}, Zhou et al., 2018 \cite{zhou2018nonconvex}, Lihua et al., 2017 \cite{lei2017nonconvex}, Nguyen et al., 2017 \cite{nguyen2017MBSARAH}, etc), and methods for composite setting remain limited (Wang et al., 2019 \cite{wang2019spiderboost}, Reddi et al., 2016b \cite{j2016nonconvex}). Moreover, most stochastic methods have been developed for convex optimization, extending them to nonconvex problems may require rigorous conditions or be limited to specific cases (see, e.g., Zhou and Gu, 2019 \cite{zhou}, Nguyen et al., 2017 \cite{nguyen2017MBSARAH}). Nowadays, there has been a growing trend in developing new algorithmic frameworks for efficient composite nonconvex optimization.

Conventional CG approaches exhibit impressive efficiency in large-scale optimization, the SCG methods inherit pivotal properties from CG and have shown promising potential in stochastic optimization. Indeed, search directions in SCG may lose the conjugacy among non-overlap samples and we have opportunities for momentum speed-up techniques. In nonconvex case, SAGA and SVRG lead to an overall SFO complexity of $O(n^{2/3}\epsilon^{-2})$ (Reddi et al., 2016 \cite{reddi2016complexity}, Allen-Zhu et al., 2016 \cite{allen2016complexity}) and SARAH yields an overall SFO complexity in the order of $\mathcal{O}(\epsilon^{-4})$ (Nguyen et al., 2017 \cite{nguyen2017MBSARAH}). Based on these estimators, variance reduced SCG approaches (e.g. Xue et al., 2021 \cite{xue2021CG4}, Kou and Yang, 2022 \cite{kou2022cgsaga}, Yang, 2023 \cite{yang2023PBCG}) have demonstrated remarkable robustness in handling noisy estimates. Motivated by this, we aim to develop novel accelerated stochastic conjugate frameworks for solving the composite problem (\ref{problem}).

Restart leads to $d$-step quadratic convergence (Nocedal and J. Wright, 2016 \cite{jorge2006numerical}) in CG, i.e., $\|w_{k+d}-w^*\|=O\left(\|w_k-w^*\|^2\right)$. It periodically refreshes the algorithms by taking steepest descent steps. A restart should be conducted whenever two consecutive gradients are far from orthogonal, as can be measured by $|\langle \nabla f_k, \nabla f_{k-1}\rangle|/\|\nabla f_k\|^2\geq\nu$, where a common choice for $\nu$ is $0.1$. In stochastic settings, SCG carried over the restart by employing related formulas (e.g., PR, PR+, FRPR) (Jin et al., 2018 \cite{jin2018cg-svrg}, Kou and Yang, 2022 \cite{kou2022cgsaga}, Yang, 2022 \cite{yang2022cgsarahhd}) or by manually setting restart conditions (e.g., Xue et al., 2021 \cite{xue2021CG4}). In fact, these algorithms renew the iterative procedures by conducting the variance reduced steps, which relied on stochastic estimates. To speed up variance reduction, we expect to perform a new restart scheme that uses exact estimates
(the full gradients) in particular refreshing steps. Notably, this design will lead to two types of restarts: one triggered by the stochastic conjugate formulas with automation into the variance reduced steps; the other from our scheduled restart nodes in fixed frequency performing the deterministic steps.

We also attempt to refine some relevant previous work. In CG, the search direction may fail to be a descent direction. One may encounter a similar issue in SCG, as observed in the experiment results of CG-SARAH -SO and CG-SARAH+ -SO algorithms (Yang, 2022 \cite{yang2022cgsarahso}), where the ascent steps have shown up in several early epochs. In terms of theory, Yang, 2022 \cite{yang2022cgsarahso} and Yang, 2022 \cite{yang2022cgsarahhd} omitted several positive terms with respect to $\mathbb{E}[\|d_k\|^2]$ in the proof procedure, which can be insufficient for deductions (Page 10 and Page 12, respectively). To improve associated demerit, we introduce the stochastic version of strong Wolfe conditions in the paper. 

Due to the computational burden in the line search, we refer to a class of probabilistic switching methods (e.g., Li et al., 2021 \cite{li2021switchpage}) and proposed a practical variant with rigorous complexity analysis. Similar to Nguyen et al. (2017) \cite{nguyen2017MBSARAH} and Tran-Dinh et al. (2022) \cite{tran2022hybrid}, our analysis has been extended to gradient dominant cases (Wang et al, 2019 \cite{wang2019spiderboost}). At last, we emphasize that our frameworks can be further extended to the preconditioned SCG method (e.g., Li et al., 2018 \cite{li2018CG3}).

The main contributions are summarized below, while the key complexity results are listed in Table \ref{table1}.
\begin{enumerate}
	\item[(a)] We propose a new stochastic conjugate framework with momentum acceleration technique to solve nonconvex composite problems. Our algorithm relies on mini-batch SARAH estimators and remains a variance reduction method. Different from general stochastic conjugate algorithms, we impose additional constraints on the update rules to mitigate possible adverse effects from outliers.
	\item[(b)] We develop a new restart scheme and present the second stochastic conjugate framework. It periodically refreshes the updates at a fix frequency. Theoretical analysis has proven its potential to outperform our first one.
	\item[(c)] We prove that both algorithms converge linearly with identical rate constants but distinct (convergence) radii. We extend our analysis to the gradient dominant case. Also, we discuss the hyper-parameters involved in the line search.
	\item[(d)] We propose a practical variant in order to reduce the computational burden associated with line search. We also provide convergence analysis for this variant.
\end{enumerate}

\begin{table}
	\setlength{\abovecaptionskip}{0pt}
	\setlength{\belowcaptionskip}{1pt}
	\renewcommand{\arraystretch}{1.3} 
	\centering
	\vspace{4pt}
	\resizebox{0.9\textwidth}{!}{
    \begin{tabular}{lccc}
		\hline
		\multicolumn{1}{c}{Algorithms}&\multicolumn{1}{c}{\qquad SFO} &\multicolumn{1}{c}{\qquad PO} &\multicolumn{1}{c}{\qquad Composite}
		\\ \hline\noalign{\smallskip}
		SCSG \cite{lei2017nonconvex} \cite{ghadimi2016batch} & \qquad $\mathcal{O}(n+n^{2/3}\epsilon^{-2})$ &  \qquad NA & \qquad \ding{55}\\
		SNVRG \cite{zhou2020SNVRG} \cite{ghadimi2016batch} & \qquad $\mathcal{O}\left((n+n^{1/2}\varepsilon^{-2})log(n)\right)$ &  \qquad NA & \qquad \ding{55}\\
		ProxGD \cite{ghadimi2016batch} & \qquad ${\mathcal{O}}(n\epsilon^{-2})$ &  \qquad ${\mathcal{O}}(\epsilon^{-2})$ & \qquad $\checkmark$\\
		ProxSVRG+ \cite{li2018proxsvrg+} & \qquad $\mathcal{O}(n+n^{2/3}\epsilon^{-2})$ & \qquad ${\mathcal{O}}(\epsilon^{-2})$ & \qquad $\checkmark$ \\
		ProxSVRG \cite{j2016nonconvex} & \qquad $\mathcal{O}(n+n^{2/3}\epsilon^{-2})$ & \qquad ${\mathcal{O}}(\epsilon^{-2})$ & \qquad$\checkmark$ \\
		ProxSAGA \cite{j2016nonconvex} & \qquad $\mathcal{O}(n+n^{2/3}\epsilon^{-2})$ & \qquad ${\mathcal{O}}(\epsilon^{-2})$ & \qquad$\checkmark$ \\
		ProxSARAH \cite{phamproxsarah} & \qquad $\mathcal{O}\left(n+n^{1/2}\epsilon^{-2}\right)$ & \qquad ${\mathcal{O}}(\epsilon^{-2})$ & \qquad $\checkmark$ \\
		SPIDER \cite{fang2018spider} & \qquad $\mathcal{O}\left(n+n^{1/2}\epsilon^{-2}\right)$ & \qquad NA & \qquad \ding{55} \\
		Prox-SpiderBoost \cite{wang2019spiderboost} & \qquad $\mathcal{O}\left(n+n^{1/2}\epsilon^{-2}\right)$ & \qquad ${\mathcal{O}}(\epsilon^{-2})$ & \qquad$\checkmark$\\
		ProxHSGD \cite{tran2022hybrid} & \qquad $\mathcal{O}\left(n+\varepsilon^{-3}\right)$ & \qquad ${\mathcal{O}}(\epsilon^{-3})$ & \qquad$\checkmark$\\
		\textcolor{blue}{This paper} & \qquad \textcolor{blue}{ $\mathcal{O}\left((n+\frac{\hat{\beta}^4}{(1-\hat{\beta})^4}\epsilon^{-2})\log(\frac{1}{\epsilon})\right)$} & \qquad \textcolor{blue}{ $\mathcal{O}\left(\frac{\hat{\beta}^2}{(1-\hat{\beta})^2}\epsilon^{-1}\log(\frac{1}{\epsilon})\right)$} & \qquad $\checkmark$\\
		\hline
	\end{tabular}}
	\vspace{4pt}
	\caption{A comparison of stochastic first-order oracle complexity (SFO) and proximal oracle complexity (PO) for nonsmooth nonconvex optimization: the results marked in blue are the contributions of the paper. SFO and PO are measured by the number of oracle calls (gradient evaluations and proximal projections) required to achieve an $\epsilon$-optimal solution. Note that all the complexity bounds here must rely on the Lipschitz constant $L$ of the smooth components and the initial functional gap $P(w_0)-P(w_{\star})$. For the sake of convenience, we omit such related quantities in the complexity bounds.}
	\label{table1}
\end{table}

\section{Preliminaries}
In the paper, we focus on stochastic estimators in the mini-batch setting. To better learn about the preliminaries, we first offer a review of biased SARAH estimators in Section \ref{SARAH} and a description of stochastic conjugate gradients in Section \ref{Conjugate}. Then, we provide the fundamental assumptions of the paper in Section \ref{Assumption}.

\subsection{Biased SARAH Estimators and Properties}\label{SARAH}
SARAH was first proposed by Nguyen et al. (2017) \cite{nguyen2017SARAH} for dealing with convex finite-sum problems. It addresses the issue remained in SVRG (Johnson et al., 2013 \cite{johnson2013SVRG}) that, iterates may have increasing trend or oscillating trend for each inner loop.

Indeed, SARAH can be viewed as a variant of SVRG. To be specific, a mini-batch SARAH estimator of gradient $\nabla f$ at $w_k^{(s)}$ ($k\geq1$) can be computed as follows:
\begin{equation}
\label{v}
v_{k}^{(s)}=\nabla f_{\mathcal{B}_k}(w_k^{(s)})-\nabla f_{\mathcal{B}_k}(w_{k-1}^{(s)})+v_{k-1}^{(s)},
\end{equation}
where $\nabla f_{\mathcal{B}_k}(w_k^{(s)})=\frac{1}{b}\sum_{i\in \mathcal{B}_k}\nabla f_i(w_k^{(s)}),\nabla f_{\mathcal{B}_k}(w_{k-1}^{(s)})=\frac{1}{b}\sum_{i\in \mathcal{B}_k}\nabla f_i(w_{k-1}^{(s)})$ and $\mathcal{B}_k\subset[n]$ with cardinality $|\mathcal{B}_k|=b$. Clearly, each evaluation of $v_{k}^{(s)}$($k\geq1$) requires two gradient evaluations of $\nabla f_{\mathcal{B}_k}$ at $w_k^{(s)}$ and $w_{k-1}^{(s)}$, which corresponds to $2b$ units of workload. Though the mini-batch SARAH estimator $v_{k}^{(s)}$ (\ref{v}) is biased, it still possesses a key property of variance reduction (see in Nguyen et al., 2017 \cite{nguyen2017SARAH}).

Next, we review some basic properties of the estimators for solving (\ref{problem}). The following Lemma \ref{lemma1} was proved in Pham et al. (2020) \cite{phamproxsarah}.
\begin{lemma}
	\label{lemma1}
	Assume that $v_{k}^{(s)}$ is updated  by (\ref{v}). Then, for any $s\geq0$, $k\geq1$,
	\begin{equation}
	\label{lemma1_1}
	\begin{aligned}
	\mathbb{E}\left[\|v_{k}^{(s)}-v_{k-1}^{(s)}\|^{2}\mid\mathcal{F}_{k}\right] &=\frac{n\left(b-1\right)}{b\left(n-1\right)}\cdot\|\nabla f(w_k^{(s)})-\nabla f(w_{k-1}^{(s)})\|^{2}\\
	&+\frac{n-b}{b\left(n-1\right)}\cdot\frac{1}{n}\sum_{i=1}^{n}\|\nabla f_{i}(w_k^{(s)})-\nabla f_{i}(w_{k-1}^{(s)})\|^{2},
	\end{aligned}
	\end{equation}
	where $\mathcal{F}_k=\sigma(w_0^{(s)},w_1^{(s)},\cdots,w_k^{(s)})$ is the $\sigma$-algebra and $\mathcal{F}_0=\mathcal{F}_1=\sigma(w_0^{(s)})$.
\end{lemma}

\vspace{4pt}
If we choose a full batch of $b=n$, we indeed use an exact gradient estimate and the second term in (\ref{lemma1_1}) vanishes, SARAH reduces to the GD algorithm.

Let us emphasize that, in the paper, when $k=0$, we choose $v_0^{(s)}=\nabla f(w_0^{(s)})$. Alternatively, one can choose $v_0^{(s)}=\frac{1}{b_2}\sum_{j\in \hat{\mathcal{B}}}\nabla f_j(w_0^{(s)})$ to utilize the inexact SARAH (i.e. ISARAH) estimators, where $\hat{\mathcal{B}}$ denotes another random mini-batch, independent of $\mathcal{B}_k$ ($k\geq1$) with $|\hat{\mathcal{B}}|=\hat{b}$. For detailed information of the ISARAH estimators and the related properties, we refer you to see Nguyen et al. (2021) \cite{nguyen2021ISARAH}.

\subsection{Stochastic Conjugate Gradient}\label{Conjugate}
For clarity, we briefly discuss the iterative updates of conventional conjugate gradient (CG) method (Nocedal and J. Wright, 2006 \cite{jorge2006numerical}) in the following steps:
\begin{align}
w_{k}&=w_{k-1}+\eta_{k-1}d_{k-1},\\
d_{k}&=-\nabla f(w_{k})+\beta_{k}d_{k-1}\label{deterministic_d},
\end{align}
where $\beta_{k}$ is conjugate parameter, and one may set $d_0=\nabla f(w_0)$ in practice. By substituting the exact $\nabla f$ with its stochastic estimate, we derive stochastic conjugate gradient, i.e., 
\begin{equation}
\label{d}
d_{k}^{(s)}=-g(w_{k}^{(s)}, \nu_k)+\beta_{k}d_{k-1}^{(s)},
\end{equation}
where $g(w_{k}^{(s)}, \nu_k)$ is the stochastic estimate at $w_k^{(s)}$.

In vanilla SCG method, one set $g(w_{k}^{(s)}, \nu_k)=\nabla f_{i_k}(w_{k}^{(s)})$ with $i_k$ selected randomly from $\{1, \ldots, n\}$. Due to the significant variance, the step size $\eta_k \propto 1 / \sqrt{k}$ may be forced to employed for sublinear convergence rate. Hence, it's adequate to use variance reduced estimates in (\ref{d}). Based on SVRG estimator, Jin et al. (2018) \cite{jin2018cg-svrg} and Xue et al. (2021) \cite{xue2021CG4} choose for $g(w_{k}^{(s)}, \nu_k)=\nabla f_{S}(w_k^{(s)})-\nabla f_S(\tilde{w})+v_0$, where $\tilde{w}$ is the snapshot. Kou and Yang (2022) \cite{kou2022cgsaga} select it relying on SAGA estimator, that is $g(w_{k}^{(s)}, \nu_k)=\nabla f_{S}(w_k)-\mu_{S}+\mu_{k-1}$, where $\mu_{S}=\frac{1}{\left|S\right|} \sum_{j \in S} \nabla f_j(w_{k-1})$. Later, Yang (2022) \cite{yang2022cgsarahso} developed another variance reduced conjugate gradient through the SARAH estimator, with $g(w_{k}^{(s)}, \nu_k)=\nabla f_{\mathcal{B}_k}(w_k^{(s)})-\nabla f_{\mathcal{B}_k}(w_{k-1}^{(s)})+v_{k-1}^{(s)}$ opted for (\ref{d}).

Prevalent options for conjugate parameter $\beta_{k}$ include but not limit to Fletcher-Reeves (FR), Polak-Ribiere (PR), Hestenes-Stiefel (HS) and Dai-Yuan (DY) formulas. Take Hestenes-Stiefel (HS) formula as an introductory example, which is not commonly used in stochastic settings, the $\beta_{k}$ in (\ref{d}) is computed by
$$
\beta_{k}^{HS}=\frac{\langle v_{k}^{(s)}, v_{k}^{(s)}-v_{k-1}^{(s)}\rangle}{\langle d_{k-1}^{(s)}, v_{k}^{(s)}-v_{k-1}^{(s)}\rangle}.
$$

In the paper, we mainly focus on the variance reduction stochastic method and FR, FR-PR update formulas, which will be specified in Section \ref{Algorithm}.

\subsection{Fundamental Assumptions}\label{Assumption}
Now, we provide some fundamental assumptions in the paper.\\
\textbf{Assumption 1} ({\bf $L$-average smoothness}). $f$ is $L$-smooth on average over any compact set of its domain, i.e., there exists $L>0$, for all $w, w^{\prime} \in \operatorname{dom}(f)$
\begin{equation}
\label{L3}
\frac{1}{n} \sum_{i=1}^n\left\|\nabla f_i(w)-\nabla f_i(w^{\prime})\right\|^2 \leq L^2\|w-w^{\prime}\|^2.
\end{equation}
$L$-average smoothness (\ref{L3}) is weaker than the individual $L$-smoothness on each $f_i$. The latter can imply the former, but not vice versa. Individual $L_i$-smoothness can achieve (\ref{L3}) as well when $\frac{1}{n} \sum_{i=1}^n L_i^2\leq L^2$. Employing (\ref{L3}), it's sufficient to derive 
\begin{equation}
\label{L}
\left\|\nabla f(w)-\nabla f(w^{\prime})\right\|^2 \leq L^2\|w-w^{\prime}\|^2.
\end{equation}
 Note that $L$-smoothness of $f$ leads to a well-known bound
\begin{equation}
\label{L2}
f(w^{\prime}) \leq f(w)+\langle\nabla f(w), w^{\prime}-w\rangle+\frac{L}{2}\|w^{\prime}-w\|^2, \quad \forall w, w^{\prime} \in \operatorname{dom}(f).
\end{equation}

\vspace{4pt}
\noindent
\textbf{Assumption 2} ({\bf Bounded conjugate parameter}). There exists a uniform upper bound, $\hat{\beta}<1$, such that
\begin{equation}
\label{beta}
\beta_{k+1}=\frac{\|v_{k+1}^{(s)}\|^2}{\|v_{k}^{(s)}\|^2} \leq \hat{\beta}, \quad \forall k, s \geq 0.
\end{equation}
Assumption 2 is commonly employed in stochastic conjugate algorithms (e.g. Xue et al., 2021 \cite{xue2021CG4}, Kou and Yang, 2022 \cite{kou2022cgsaga}, Yang, 2023 \cite{yang2023PBCG}), it ensures an uniform upper bound for the conjugate parameter of FR updates.

In Algorithm \ref{alg1}, we additionally need the following bounded deviation condition.
\noindent
\textbf{Assumption 3} ({\bf Bounded deviation}). The estimators $v_{k}^{(s)}$ in (\ref{v}) has a bounded deviation from $\nabla f$, e.g., there exists a uniform $\sigma>0$, such that  
\begin{equation}
\label{deviation}
\mathbb{E}\left[\|v_{k}^{(s)}-\nabla f(w_k^{(s)})\|^2\right] \leq \sigma^2, \quad \forall k, s \geq 0
\end{equation}

The bounded variance assumption is standard in the analysis of non-convex optimization problems (e.g. Wang et al., 2017 \cite{wang2017deviation}, Zhou et al. 2018 \cite{zhou2018deviation2}). According to Lemma 2 in Nguyen et al. (2017) \cite{nguyen2017MBSARAH} and Lemma 1 in Pham et al. (2020) \cite{phamproxsarah}, we can have a parallel and reliable assumption in terms of the bounded deviation.

\section{Algorithms}\label{Algorithm}
In this section, we introduce two new stochastic conjugate frameworks for solving problem (\ref{problem}). The first one is described in Algorithm \ref{alg1} and abbreviated by Acc-Prox-CG-SARAH algorithm. The second adopts a novel restart scheme and is referred to as Acc-Prox-CG-SARAH-RS algorithm, as outlined in Algorithm \ref{alg2}.

\begin{algorithm}[h]
	\caption{Acc-Prox-CG-SARAH}\label{alg1}
	\begin{algorithmic}[1]
		\Require initial point $\widetilde{w}_0$, epoch length $m$, mini-batch size $b$, momentum parameter $\gamma \in(0,1]$, threshold parameters $\eta_2>0$, $\beta_o>0$, $\rho>0$.
		\State $h^{(0)}=\nabla f(\widetilde{w}_0)$
		\For{$s = 1,2,...,S$}
		\State $w_0^{(s)}=\widetilde{w}_{s-1}$
		\State $v_0^{(s)}=\frac{1}{n} \sum_{i \in N} \nabla f_i(w_0^{(s)})=\nabla f(w_0^{(s)})$
		\State $d_0^{(s)}=-h^{(s-1)}$
		\State Call the line search and find $\widetilde{\eta}$ satisfying (\ref{Wolfe1}) and (\ref{Wolfe2})
		\State Set $\eta=\min\left\{\widetilde{\eta}, \eta_2 \right\}$
		\State \label{Step5} $y_0^{(s)}=\operatorname{prox}_{\eta \varphi}(w_0^{(s)} + \eta d_0^{(s)})$
		\State \label{Step1} $w_1^{(s)} =\left(1-\gamma\right) w_0^{(s)}+\gamma y_0^{(s)}$
		\For{$k = 1,2,...,m-1$}
		\State Pick mini-batch $\mathcal{B}_k \subset\{1, \ldots, n\}$ of size $b$ uniformly at random
		\State Update mini-batch SARAH estimator $v_{k}^{(s)}$ by
		\State $$v_{k}^{(s)}=\nabla f_{\mathcal{B}_k}(w_k^{(s)})-\nabla f_{\mathcal{B}_k}(w_{k-1}^{(s)})+v_{k-1}^{(s)}$$
		\State Calculate $\beta_k$ by AFR (\ref{AFR}) or FRPR (\ref{FRPR}) formula
		\State Compute the search direction $d_k^{(s)}$ by
		\State $$d_k^{(s)}=-v_{k}^{(s)}+\beta_kd_{k-1}^{(s)}$$
		\State Call the line search and find $\widetilde{\eta}$ satisfying (\ref{Wolfe1}) and (\ref{Wolfe2})
		\State Set $\eta=\min\left\{\widetilde{\eta}, \eta_2 \right\}$
		\State Update iterates by
		\State \label{Step6} $y_{k}^{(s)}=\operatorname{prox}_{\eta \varphi}(w_k^{(s)} + \eta d_k^{(s)})$
		\State \label{Step2} $w_{k+1}^{(s)} =\left(1-\gamma\right) w_k^{(s)}+\gamma y_{k}^{(s)}$
		\EndFor
		\State $h^{(s)}=v_m^{(s)}$
		\State $\widetilde{w}_s=w_m^{(s)}$
		\EndFor
	\end{algorithmic}
\end{algorithm}

\begin{algorithm}[h]
	\caption{Acc-Prox-CG-SARAH-RS}\label{alg2}
	\begin{algorithmic}[1]
		\Require initial point $\widetilde{w}_0$, epoch length $m$, mini-batch size $b$, momentum parameter $\gamma \in(0,1]$, threshold parameters $\eta_2>0$, $\beta_o>0$, $\rho>0$.
		\For{$s = 1,2,...,S$}
		\State $w_0^{(s)}=\widetilde{w}_{s-1}$
		\State $v_0^{(s)}=\frac{1}{n} \sum_{i \in N} \nabla f_i(w_0^{(s)})=\nabla f(w_0^{(s)})$
		\State \textcolor{blue}{$d_0^{(s)}=-v_0^{(s)}=-\nabla f(w_0^{(s)})$}
		\State Call the line search and find $\widetilde{\eta}$ satisfying (\ref{Wolfe1}) and (\ref{Wolfe2})
		\State Set $\eta=\min\left\{\widetilde{\eta}, \eta_2 \right\}$
		\State \label{Step7} $y_0^{(s)}=\operatorname{prox}_{\eta \varphi}(w_0^{(s)} + \eta d_0^{(s)})$
		\State \label{Step3}$w_1^{(s)} =\left(1-\gamma\right) w_0^{(s)}+\gamma y_0^{(s)}$
		\For{$k = 1,2,...,m-1$}
		\State Pick mini-batch $\mathcal{B}_k \subset\{1, \ldots, n\}$ of size $b$ uniformly at random
		\State Update mini-batch SARAH estimator $v_{k}^{(s)}$ by
		\State $$v_{k}^{(s)}=\nabla f_{\mathcal{B}_k}(w_k^{(s)})-\nabla f_{\mathcal{B}_k}(w_{k-1}^{(s)})+v_{k-1}^{(s)}$$
		\State Calculate $\beta_k$ by AFR (\ref{AFR}) or FRPR (\ref{FRPR}) formula
		\State Compute the search direction $d_k^{(s)}$ by
		\State $$d_k^{(s)}=-v_{k}^{(s)}+\beta_kd_{k-1}^{(s)}$$
		\State Call the line search and find $\widetilde{\eta}$ satisfying (\ref{Wolfe1}) and (\ref{Wolfe2})
		\State Set $\eta=\min\left\{\widetilde{\eta}, \eta_2 \right\}$
		\State Update iterates by
		\State \label{Step8} $y_{k}^{(s)}=\operatorname{prox}_{\eta \varphi}(w_k^{(s)} + \eta d_k^{(s)})$
		\State \label{Step4}$w_{k+1}^{(s)} =\left(1-\gamma\right) w_k^{(s)}+\gamma y_{k}^{(s)}$
		\EndFor
		\State $\widetilde{w}_s=w_m^{(s)}$
		\EndFor
	\end{algorithmic}
\end{algorithm}

Both frameworks remain variance reduction methods. We choose mini-batch settings of SARAH estimators (Nguyen et al., 2017 \cite{nguyen2017MBSARAH}) to better use the parallel processing power. The batch size $b$ is chosen exponentially proportional to $n$. The conjugate parameter $\beta$ is computed by adaptive FR (AFR) formula or FR-PR formula. Through FR formula, we compute the parameter $\beta$ at $w_k^{(s)}$ by
\begin{equation}
\label{FR}
\beta_{k}^{FR} =\frac{\|v_{k}^{(s)}\|^{2}}{\|v_{k-1}^{(s)}\|^{2}}. 
\end{equation}
Though the convergence of FR method has been rigorously established, their numerical results have potential inefficiencies. To control $\beta$ adaptively, we multiply another coefficient $\rho$ to $\beta_{k}^{FR}$ (\ref{FR}) and bound its maxima by a constant $\beta_o$. Then, we have our adaptive FR (AFR) update formula as
\begin{equation}
\label{AFR}
\beta_k^{AFR}=\min \left\{\beta_o, \rho\beta_{k}^{FR}\right\}.
\end{equation}
Empirically, we set $\rho\in\{0.8, 1\}$. As a pivotal variant of FR, the PR method exhibits strong robustness and inspiring efficiency in applications, which is derived by
$$
\beta_{k}^{PR}=\frac{\langle v_{k}^{(s)}, v_{k}^{(s)}-v_{k-1}^{(s)}\rangle}{\|v_{k-1}^{(s)}\|^2}. 
$$
In terms of theory, the PR method may not converge. To leverage its efficiency, we opt for a hybrid formula of FR-PR for computing $\beta$ at $w_k^{(s)}$, i.e.,
\begin{equation}
\label{FRPR}
\beta_k^{FRPR}=\left\{\begin{array}{cc}
-\beta_k^{F R}, & \text { if } \beta_k^{P R}<-\beta_k^{F R}, \\
\beta_k^{P R}, & \text { if }\left|\beta_k^{P R}\right| \leq \beta_k^{F R}, \\
\beta_k^{F R}, & \text { if } \beta_k^{P R}>\beta_k^{F R}.
\end{array}\right.
\end{equation}
From (\ref{FRPR}), it's sufficient to observe that $|\beta_{k}^{FRPR}|\leq\beta_{k}^{\mathrm{FR}}$.

In gradient steps of $f$, we perform the strong Wolfe condition in a stochastic manner (Jin et al., 2018 \cite{jin2018cg-svrg}, Kou and Yang, 2022 \cite{kou2022cgsaga}) to preserve the (batch) descent property, as well to search for appropriate step sizes, i.e., 
\begin{align}
f_{\mathcal{B}_k}(w_{k}^{(s)}+\widetilde{\eta} d_{k}^{(s)}) &\leq f_{\mathcal{B}_k}(w_{k}^{(s)})+c_{1}\widetilde{\eta}\langle\nabla f_{\mathcal{B}_k}(w_{k}^{(s)}),d_{k}^{(s)}\rangle,  \label{Wolfe1}\\
\left|\langle v_{k+1}^{(s)}, d_{k}^{(s)}\rangle\right| &\leq-c_{2}\langle v_{k}^{(s)}, d_{k}^{(s)}\rangle,  \label{Wolfe2}
\end{align}
where $0<c_1<c_2<1$. The searched step size $\widetilde{\eta}$ satisfies conditions (\ref{Wolfe1}) and (\ref{Wolfe2}). To mitigate the effects of some outliers in samples, we add a constraint on step sizes, where we bound $\eta=\min\left\{\widetilde{\eta}, \eta_2 \right\}$ ($\eta_2>0$ is predefined).

In proximal steps, we proceed to choose $\eta$ which leads to the proximal stochastic gradient descent (PSGD) at Step \ref{Step5}, \ref{Step6}, \ref{Step7}, \ref{Step8} (Ghadimi and Lan, 2013 \cite{ghadimi2013proximal}), i.e., 
$$y_{k}^{(s)}=\operatorname{prox}_{\eta \varphi}(w_k^{(s)} + \eta d_k^{(s)}).$$

The Nesterov's momentum (Nesterov, 1983 \cite{nesterov1983acc}) at Step \ref{Step1}, \ref{Step2}, \ref{Step3}, \ref{Step4} is applied to accelerate our first-order optimization with weight $\gamma \in(0,1]$. The first update on $y_k^{(s)}$ is a standard PSGD step, while the second one on $w_k^{(s)}$ is a momentum step. The mini-batch SARAH estimator $v_{k}^{(s)}$ is evaluated at the averaged point $w_k^{(s)}$, instead of the middle iterate $y_k^{(s)}$. The convex combination of $w_k^{(s)}$ and $y_k^{(s)}$ with momentum coefficient $\gamma$ yields $w_{k+1}^{(s)}$, which can be re-expressed as
$$
w_{k+1}^{(s)}=(1-\gamma)w_{k}^{(s)}+\gamma\text{prox}_{\eta\varphi}(w_{k}^{(s)}+\eta d_{k}^{(s)}).
$$
Equivalently, it's an exponentially weighted average of middle iterates $\{y_k^{(s)}\}$ to $w_{k}^{(s)}$,
$$
w_{k}^{(s)}=(1-\gamma)^{k}w_{0}^{(s)}+\gamma\sum_{i=1}^{k}(1-\gamma)^{i-1}y_{k-i}^{(s)}.
$$
If $\varphi=0$ and we set $\gamma=1$, i.e., for non-composite optimization without momentum acceleration, Algorithm \ref{alg1} is reduced to a simplified variant that shares similarities with those presented in Yang, 2022 \cite{yang2022cgsarahso} and Yang, 2022 \cite{yang2022cgsarahhd}.

Note that  $H_{\eta}(w_k^{(s)})=\frac{1}{\eta}(w_k^{(s)}-\operatorname{prox}_{\eta \varphi}(w_k^{(s)}+\eta d_k^{(s)}))$ is a stochastic approximation to the gradient mapping at $w_k^{(s)}$, i.e.,  $\mathcal{G}_{\eta}(w_k^{(s)})=\frac{1}{\eta}(w_k^{(s)}-\operatorname{prox}_{\eta \varphi}(w_k^{(s)}-\eta \nabla  f(w_k^{(s)})))$. Therefore, it's sufficient to rewrite the key steps as
$$
w_{k+1}^{(s)}=w_k^{(s)}-\eta \gamma H_{\eta}(w_k^{(s)}).
$$
If $\varphi=0$, then $H_{\eta}(w_k^{(s)})= v_{k}^{(s)}$, which reduces to the SARAH approximation to $\nabla f$. The product $\eta \gamma$ is considered as an effective step size or a combined step size.

The FR-PR formula (\ref{FRPR}) automatically restarts the optimization with $v_{k}^{(s)}$ for the inner loops, whenever $v_{k}^{(s)}\approx v_{k-1}^{(s)}$ and $\left|\beta_k^{P R}\right| \leq \beta_k^{FR}$. Apart from that, it's also encouraging to impose constraints on $\beta_k^{FRPR}$ (\ref{FRPR}), e.g., we can set $\beta_k=0$ whenever $\beta_k^{FRPR}$ (\ref{FRPR}) exceeds certain predefined threshold values.

In Algorithm \ref{alg2}, we additionally apply a novel restart scheme of enforcing $d_0^{(s)}=\nabla f(w_0^{(s)})$ for all steps $s\geq0$. Then, the epoch length $m>0$ corresponds to a fixed restart interval. Though $v_{k}^{(s)}$ (\ref{v}) replaces the role of $\nabla f$ in stochastic settings, Algorithm \ref{alg2} in fact uses $\nabla f$ to renew the search directions at beginning of each epoch, instead of the estimate $v_{k}^{(s)}$. It indicates that we actually possess two opportunities to refresh the algorithm: one through FR-PR formula, which automatically refreshes the poor search directions of stochastic conjugate gradients; the other via our restart scheme in fixed frequency, regularly mitigating the variance (deviation) accumulated from stochastic optimization (the mechanism is theoretically specified in Section \ref{Convergence}). We call this novel restart scheme as deterministic restart.

\section{Convergence Analysis}\label{Convergence}
\label{section}
We analyze the convergence result with any $\beta_{k}$ satisfying $|\beta_k|\leq\beta_k^{FR}$ (For reference, we recommend you to see Lemma 3.1 in Gilbert et al., 1992 \cite{gilbert1992proof}). Let us emphasize that the convergence analysis in Algorithm \ref{alg2} dose not require Assumption 3.
\subsection{Optimality Conditions}
Since $\varphi$ is proper, closed and convex, its proximal mapping $\operatorname{prox}_{\eta \varphi}(\cdot)$ is well defined, i.e., for any $w \in \operatorname{dom}(f) \cap \operatorname{dom}(\varphi)\neq\emptyset$, the proximal mapping $\operatorname{prox}_{\eta \varphi}(w)$ exists and is unique (Nesterov, 2003 \cite{nesterov2003introductory}). If $\varphi=0$, it becomes the identity mapping.

Assume that $P$ is bounded from below. When $f$ is nonconvex, the first order (stationary) optimality condition of problem (\ref{problem}) can be stated as 
\begin{equation}
\label{optimality}
0 \in \partial P(w_{\star}) \equiv \nabla f(w_{\star})+\partial \varphi(w_{\star}),
\end{equation}
where $\partial \varphi$, $\partial P$ denote corresponding subdifferentials at $w_{\star}$, and $w_{\star}$ is called a stationary point of $P$. In order to analyze the convergence results for nonsmooth nonconvex problems, we define the gradient mapping of $P$ at $w$ as

\begin{equation}
\label{mapping}
\mathcal{G}_\eta(w)=\frac{1}{\eta}\left(w-\operatorname{prox}_{\eta \varphi}\left(w-\eta \nabla f(w)\right)\right), \quad \forall \eta>0.
\end{equation}
If $\varphi=0$, this mapping reduces to $\mathcal{G}_\eta(w)=\nabla f(w)$. By using $\mathcal{G}_\eta$ (\ref{mapping}) as the convergence criterion, the (stationary) optimality condition (\ref{optimality}) is naturally specified as
$$
\left\|\mathcal{G}_\eta(w_{\star})\right\|^2=0.
$$

In stochastic settings, it's common practice to bound the total number of iterations, $\mathcal{T}$, to achieve an $\epsilon$-optimal solution to $w_{\star}$ ($\epsilon$ is predefined). For clarity, we define the $\epsilon$-optimal solution of (\ref{problem}) as follows:
\vspace{4pt}
\begin{definition}
	Given predefined accuracy $\epsilon>0$, the iterate $w_{\mathcal{T}}$ is said to be an $\epsilon$-optimal solution of (\ref{problem}), if
	\begin{equation}
	\label{optimality2}
	\mathbb{E}\left[\|\mathcal{G}_\eta(w_\mathcal{T})\|^2\right]\leq\epsilon.
	\end{equation}
\end{definition}
\vspace{4pt}
In the paper, the complexity analysis is established based on condition (\ref{optimality2}) for Incremental First-Order Oracle Model (Agarwal and Bottou, 2014 \cite{agarwal}).

\subsection{A Useful Lemma}
To begin with, %we define an auxiliary function as follows:
%\begin{definition}
	let us define an auxiliary function
	$$\Phi(x)=(1+x)/(1-x).$$ 
It is obvious that $\Phi(x)$ is a monotone increasing function. 
 
%\end{definition}
We present an important result that is essential in the following convergence analysis. It demonstrates how the line search (\ref{Wolfe2}) affects the stochastic search directions and update procedures. With slight modifications to Theorem 2 in Jin et al. (2018) \cite{jin2018cg-svrg}, we estimate the upper bound of the expected values of $\|d_k^{(s)}\|^2$.
\begin{lemma}
	\label{lemma2}
	Suppose that Assumption 2 and Assumption 3 hold. Let Algorithm \ref{alg1} be implemented with step size $\widetilde{\eta}$ that satisfies condition (\ref{Wolfe2}). For a chosen $c_2\in(0, 1)$, assume that $\Phi(c_2)\leq\alpha$ ($\alpha>1$). Then, for all $s, k\geq0$, there exists a constant $\tau>0$ such that
	$$\mathbb{E}\left[\|d_k^{(s)}\|^2\right]\leq g(k)\mathbb{E}\left[\|v_0^{(s)}\|^2\right]+\hat{\beta}^{2k}\tau\sigma^2,$$
	where 
	$$g(k)=\frac{\alpha}{1-\hat{\beta}}\times\hat{\beta}^{k}-\frac{\alpha-1+\hat{\beta}}{1-\hat{\beta}}\times\hat{\beta}^{2k}.$$
	In Algorithm \ref{alg2} with $d_{0}^{(s)}=v_0^{(s)}$ ($s\geq1$), for all $s, k\geq0$ we have
	$$\mathbb{E}\left[\|d_k^{(s)}\|^2\right]\leq g(k)\mathbb{E}\left[\|v_0^{(s)}\|^2\right].$$
\end{lemma}
\proof  Deferred to the Appendix A.

\subsection{Analysis of the Inner Loops}
Now, we provide the analysis of the inner loops as the foundation of our convergence results. We break the analysis into two cases and proceed step by step.
\begin{enumerate}
	\item[*] \textbf{Case 1:} Consider Algorithm \ref{alg1}, \ref{alg2} with one single outer loop where $S=1$.
	\item[*] \textbf{Case 2:} Consider Algorithm \ref{alg1}, \ref{alg2} with multiple outer loops where $S>1$.
\end{enumerate}

We start by proposing subsequent Lemma \ref{lemma3} to analyze the \textbf{Case 1}. It then will be extended to the \textbf{Case 2} and summarized in Theorem \ref{theorem1}. For brevity, we leave out the superscripts at top right of the notations.
\begin{lemma}
	\label{lemma3}
	 Let Algorithm \ref{alg1}, \ref{alg2} be implemented with $S=1$. Assume that during iterations $\eta\in\{\eta_1, \eta_2\}$. Then, for any $m>0$,
	$$
	\begin{aligned}
	&\mathbb{E}\left[P(w_{m+1})\right]\\
	&\leq \mathbb{E}\left[P(w_{0})\right]-\frac{\gamma\eta_1^{2}}{2}\sum_{k=0}^{m}\mathbb{E}\left[\|\mathcal{G}_{\eta}(w_{k})\|^{2}\right]+ \gamma\left(1+2\eta_2^{2}\right)\hat{\beta}^2\sum_{k=1}^{m}\mathbb{E}\left[\|d_{k-1}\|^2\right] \\
	&+\gamma\left(1+2\eta_2^{2}\right)\sum_{k=0}^{m}\mathbb{E}\left[\|\nabla f(w_{k})-v_{k}\|^{2}\right]-\frac{\gamma}{2}\left(\frac2{\eta_2}-L\gamma-3\right)\sum_{k=0}^{m}\mathbb{E}\left[\|y_{k}-w_{k}\|^{2}\right].
	\end{aligned}
	$$
\end{lemma}
\proof  Deferred to the Appendix B.

\vspace{8pt}
Subsequently, we consider Algorithm \ref{alg1}, \ref{alg2} with multiple outer loops ($S>1$) and with $s>1$. In this case, we actually analyze the inner loops within the specific $s$-th epoch. We establish our first convergence result in Theorem \ref{theorem1}.
\begin{theorem}
	\label{theorem1}
	Suppose that Assumptions 2 and 3 hold. Let $w_{\star}$ be the unique minimizer of $P$. Assume that during iterations $\widetilde{\eta}\in\{\eta_{1}, \eta_{2}\}$, and parameters $b$, $\gamma$ are chosen under a simple and suitable condition \footnote{This condition will be specified as \eqref{Thm1} in Appendix C.}. Then, for all $s>1$ in Algorithm \ref{alg1}, we have
	\begin{equation}
	\label{rate1}
	\mathbb{E}\left[\|\mathcal{G}_{\eta}(w_{m}^{(s)})\|^{2}\right]
	\leq\frac{\xi(1-\hat{\beta})^2(1+\hat{\beta})}{\gamma\alpha(1+2\eta_2^2)\hat{\beta}^2}\mathbb{E}\left[P(w_{0}^{(s)})-P(w_{\star})\right]+\xi\mathbb{E}\left[\|\nabla f(w_0^{(s)})\|^2\right]+C\xi,
	\end{equation}
	where parameters $\xi$ and $C$ are given by
	\begin{equation}
	\label{C}
	\begin{aligned}
	&\xi=\frac{(2+4\eta_2^2)\alpha\hat{\beta}^2}{(m+1)\eta_1^2(1-\hat{\beta})^2(1+\hat{\beta})},\\
	\quad
	&C=\left(\frac{1-\hat{\beta}}{\alpha}\tau+\frac{(1-\hat{\beta})^2(1+\hat{\beta})}{2\alpha\hat{\beta}^2}\right)\sigma^2.
	\end{aligned}
	\end{equation}
In Algorithm \ref{alg2} with $d_{0}^{(s)}=v_0^{(s)}$ ($s\geq1$), we have 
	$$
	\mathbb{E}\left[\|\mathcal{G}_{\eta}(w_{m}^{(s)})\|^{2}\right]
	\leq\frac{\xi(1-\hat{\beta})^2(1+\hat{\beta})}{\gamma\alpha(1+2\eta_2^2)\hat{\beta}^2}\mathbb{E}\left[P(w_{0}^{(s)})-P(w_{\star})\right]+\xi\mathbb{E}\left[\|\nabla f(w_0^{(s)})\|^2\right].
	$$
\end{theorem}

\proof Deferred to the Appendix C.

\vspace{8pt}
In fact, the analysis of the first inner loop in \textbf{Case 2} ($S>1$, $s=1$) is equivalent to that in \textbf{Case 1} ($S=1$), for the initial values $v_0=\nabla f(\widetilde{w}_0)$ and $v_0^{(1)}=\nabla f(\widetilde{w}_0)$ are equally required. Besides, Algorithm \ref{alg2} with the deterministic restart scheme periodically refreshes itself by $d_0^{(s)}=\nabla f(\widetilde{w}_s)$, hence, the technical results for inner analysis are always equivalent in \textbf{Case 1} and \textbf{Case 2} in Algorithm \ref{alg2}. 

Theorem \ref{theorem1} shows sublinear convergence rates of the inner loops of Algorithm \ref{alg1} and Algorithm \ref{alg2} with increasing $m$. For any given accuracy $\epsilon>0$, to ensure $\mathbb{E}[\|\mathcal{G}_{\eta}(w_{m}^{(s)})\|^{2}]\leq\epsilon$ within the inner loops, it's sufficient to choose $m=\mathcal{O}\left(\frac{\hat{\beta}^2}{\epsilon(1-\hat{\beta})^2}\right)$. Further, according to the condition (\ref{Thm1}), we may choose $b=\mathcal{O}\left(\frac{L^2\gamma^2m}{2-L\gamma\eta_2-3}\right)$. Then, the complexity of Algorithm \ref{alg1} and Algorithm \ref{alg2} for achieving an $\epsilon$-optimal solution within the inner loops is in the order of $n+2mb=\mathcal{O}\left(n+\frac{L^2\gamma^2\hat{\beta}^4}{\epsilon^2(2-L\gamma\eta_2-3)(1-\hat{\beta})^4}\right)$. Note that the term $C$ in (\ref{rate1}), which corresponds to a stochastic also historical deviation, does not affect the order of the convergence rate and the complexity bound. 

\paragraph{Further Discussion on $c_2$: }  Since the line search (\ref{Wolfe2}) parameter $c_2\in(0, 1)$,  we actually use a general upper bound of $\Phi(c_2)$ in our analysis. Observed from previous works, we can bound the parameter $c_2$ within certain ranges to obtain specific convergence results. One popular choice for this effective range is $c_2\in(0, \frac{1}{5})$ (e.g.,  Kou and Yang, 2022 \cite{kou2022cgsaga}). In this case, we obtain the supremum of $\alpha=2$ for $\Phi$. Thus, ineq. (\ref{Ed}) in Lemma \ref{lemma2} becomes:
\begin{equation}
\mathbb{E}\left[\|d_{k}^{(s)}\|^{2}\right]
=\left(\frac{2}{1-\hat{\beta}}\hat{\beta}^{k}-\frac{1+\hat{\beta}}{1-\hat{\beta}}\hat{\beta}^{2k}\right)\mathbb{E}\left[\|v_0^{(s)}\|^2\right]+\hat{\beta}^{2k}\tau\sigma^2 
\end{equation}
Next, ineq. (\ref{Thm_3}) in Theorem \ref{theorem1} converts into:
\begin{equation}
\begin{aligned}
\sum_{k=0}^{m-1}\mathbb{E}\left[\|d_{k}^{(s)}\|^2\right]
&\leq\frac{(1-\hat{\beta}^m)^2}{(1-\hat{\beta})^2}\mathbb{E}\left[\|\nabla f(w_0^{(s)})\|^2\right]+\frac{1-\hat{\beta}^{2m}}{1-\hat{\beta}^2}\tau\sigma^2\\
&\leq\frac{1}{(1-\hat{\beta})^2}\mathbb{E}\left[\|\nabla f(w_0^{(s)})\|^2\right]+\frac{1}{1-\hat{\beta}^2}\tau\sigma^2.
\end{aligned}
\end{equation}

By deducing further, we attain the following Corollary \ref{corollary1} that illustrates the convergence results under the specified range $c_2\in(0, \frac{1}{5})$.
\begin{corollary}
	\label{corollary1}
	Let Algorithm \ref{alg1} be implemented with step size $\widetilde{\eta}$ that satisfies condition (\ref{Wolfe2}) with $c_2\in(0, \frac{1}{5})$. Assume that during iterations $\widetilde{\eta}\in\{\eta_{1}, \eta_{2}\}$ and we choose $b, \gamma$ that satisfy condition (\ref{Thm1}). Then, for all $s>1$ in Algorithm \ref{alg1}, we have
	\begin{equation}
	\mathbb{E}\left[\|\mathcal{G}_{\eta}(w_{m}^{(s)})\|^{2}\right]
	\leq\frac{\xi_2(1-\hat{\beta})^2}{\gamma(1+2\eta_2^2)\hat{\beta}^2}\mathbb{E}\left[P(w_{0}^{(s)})-P(w_{\star})\right]+\xi_2\mathbb{E}\left[\|\nabla f(w_0^{(s)})\|^2\right]+C_2\xi_2,
	\end{equation}
	where parameters $\xi_2$ and $C_2$ are given by
	$$
	%\begin{aligned}
	 \xi_2=\frac{(2+4\eta_2^2)\hat{\beta}^2}{(m+1)\eta_1^2(1-\hat{\beta})^2}, 
	\quad
	 C_2=\left(\frac{1-\hat{\beta}}{1+\hat{\beta}}\tau+\frac{(1-\hat{\beta})^2}{2\hat{\beta}^2}\right)\sigma^2.
	%\end{aligned}
	$$
In Algorithm \ref{alg2} that is developed with $d_{0}^{(s)}=v_0^{(s)}$, we have 
	$$
	\mathbb{E}\left[\|\mathcal{G}_{\eta}(w_{m}^{(s)})\|^{2}\right]
	\leq\frac{\xi_2(1-\hat{\beta})^2}{\gamma(1+2\eta_2^2)\hat{\beta}^2}\mathbb{E}\left[P(w_{0}^{(s)})-P(w_{\star})\right]+\xi_2\mathbb{E}\left[\|\nabla f(w_0^{(s)})\|^2\right].
	$$
\end{corollary}

\subsection{Convergence Analysis}
By simply applying Theorem \ref{theorem1} to each of the outer loops, we obtain the second convergence results stated in Theorem \ref{theorem2}. Both Algorithm \ref{alg1} and Algorithm \ref{alg2} achieve linear convergence rates for the nonconvex, nonsmooth composite problem (\ref{problem}).
\begin{theorem}
	\label{theorem2}
	Suppose that Lemma \ref{lemma2}, \ref{lemma3} and Theorem \ref{theorem1} hold. Let $w_{\star}$ be the unique minimizer of $P(\cdot)$. Assume that during iterations $\widetilde{\eta}\in\{\eta_{1}, \eta_{2}\}$. Define 
	$$\delta_k=\frac{(1-\hat{\beta})^2(1+\hat{\beta})}{\gamma\alpha(1+2\eta_2^2)\hat{\beta}^2}\mathbb{E}\left[P(\widetilde{w}_k)-P(w_{\star})\right], \quad k=0, 1, \cdots, s-1,$$
	and $\delta=\operatorname*{max}_{0\leq k\leq s-1}\delta_{k}$. Assume that we choose $\xi=\frac{(2+4\eta_2^2)\alpha\hat{\beta}^2}{(m+1)\eta_1^2(1-\hat{\beta})^2(1+\hat{\beta})}<1$ and $b$, $\gamma$ under condition (\ref{Thm1}). By setting $\Delta=\frac{\xi(\delta+C)}{1-\xi}$, we have for Algorithm \ref{alg1} that
	\begin{equation}
	\label{convergence1}
	\mathbb{E}\left[\|\mathcal{G}_{\eta}(\widetilde{w}_s)\|^2\right]-\Delta\leq\xi^s\left(\|\nabla f(\widetilde{w}_0)\|^2-\Delta\right).
	\end{equation}
    Therefore, $\{\mathbb{E}\left[\|\mathcal{G}_{\eta}(\widetilde{w}_s)\|^2\right]\}$ converges linearly with rate $\xi$ to a $\Delta$-ball around zero.
    
    \noindent
In Algorithm \ref{alg2} with $d_{0}^{(s)}=v_0^{(s)}$ ($s\geq1$), by setting $\bar{\Delta}=\frac{\xi\delta}{1-\xi}$, we have
	\begin{equation}
	\label{convergence2}
	\mathbb{E}\left[\|\mathcal{G}_{\eta}(\widetilde{w}_s)\|^2\right]-\bar{\Delta}\leq\xi^s\left(\|\nabla f(\widetilde{w}_0)\|^2-\bar{\Delta}\right).
	\end{equation}
	Thus, $\{\mathbb{E}\left[\|\mathcal{G}_{\eta}(\widetilde{w}_s)\|^2\right]\}$ converges linearly with a rate of $\xi$ to a $\bar{\Delta}$-ball around zero.
\end{theorem}
\proof 
Firstly, we present the proof associated with Algorithm \ref{alg1}. According to Theorem \ref{theorem1}, it's adequate to obtain
\begin{equation}
\label{Thm2_1}
\mathbb{E}\left[\|\mathcal{G}_{\eta}(w_{m}^{(s)})\|^{2}\right]
\leq\xi\left(\frac{(1-\hat{\beta})^2(1+\hat{\beta})}{\gamma\alpha(1+2\eta_2^2)\hat{\beta}^2}\mathbb{E}\left[P(w_{0}^{(s)})-P(w_{\star})\right]+\mathbb{E}\left[\|\nabla f(w_0^{(s)})\|^2\right]+C\right).
\end{equation}
For the convenience, we define
$$\delta_k=\frac{(1-\hat{\beta})^2(1+\hat{\beta})}{\gamma\alpha(1+2\eta_2^2)\hat{\beta}^2}\mathbb{E}\left[P(\widetilde{w}_k)-P(w_{\star})\right]\quad\text{and}\quad\delta=\operatorname*{max}_{0\leq k\leq s-1}\delta_{k}.$$
Note that $w_{m}^{(s)}=\widetilde{w}_s$ and $w_0^{(s)}=\widetilde{w}_{s-1}$, we derive from (\ref{Thm2_1}) that
\begin{equation}
\mathbb{E}\left[\|\mathcal{G}_{\eta}(\widetilde{w}_s)\|^{2}\right]
\leq\xi\left(\mathbb{E}\left[\|\nabla f(\widetilde{w}_{s-1})\|^2\right]+\delta+C\right).
\end{equation}
Through the definition of $\Delta=\frac{\xi(\delta+C)}{1-\xi}$, we obtain
$$
\begin{aligned}
\mathbb{E}\left[\|\mathcal{G}_{\eta}(\widetilde{w}_s)\|^2\right]-\Delta&\leq\xi\left(\mathbb{E}\left[\|\nabla f(\widetilde{w}_{s-1})\|^2\right]-\Delta\right)\\
&\leq\xi^s\left(\|\nabla f(\widetilde{w}_{0})\|^2-\Delta\right).
\end{aligned}
$$

In Algorithm \ref{alg2} with $d_{0}^{(s)}=v_0^{(s)}$, by defining $\bar{\Delta}=\frac{\xi\delta}{1-\xi}$, we have
$$
\mathbb{E}\left[\|\mathcal{G}_{\eta}(\widetilde{w}_s)\|^2\right]-\bar{\Delta}
\leq\xi^s\left(\|\nabla f(\widetilde{w}_{0})\|^2-\bar{\Delta}\right).
$$

\vspace{8pt}
Unlike the convergence results in Nguyen et al. (2017) \cite{nguyen2017SARAH}, here, the rate constant $\xi$ in (\ref{C}) decreases as the epoch length $m$ increases. By properly `stretching' the epoch length $m$, we may obtain faster convergence results within the same linear rate order. Moreover, the convergence radii $\Delta$ and $\bar{\Delta}$ decrease as $m$ increases.

Based on Theorem \ref{theorem2}, it's adequate to choose $\Delta=\mathcal{O}\left(\epsilon\right)$ or $\bar{\Delta}=\mathcal{O}\left(\epsilon\right)$. In fact, we can choose $\Delta=\bar{\Delta}=\epsilon/4$. Thus, the total SFO complexity to achieve $\epsilon$-accuracy defined in (\ref{optimality2}) is in the order of $\mathcal{O}\left((n+\frac{L^2\gamma^2\hat{\beta}^4}{\epsilon^2(2-L\gamma\eta_2-3)(1-\hat{\beta})^4})log(\frac{1}{\epsilon})\right)$. Furthermore, the total number of calls to proximal oracle (PO) does not exceed $\mathcal{O}\left(\frac{\hat{\beta}^2}{\epsilon(1-\hat{\beta})^2}log(\frac{1}{\epsilon})\right)$.

With identical $\{\gamma, m, b\}$ and assuming $\delta=\delta_0$ (maximum sub-optimality achieved at initial $\widetilde{w}_0$), the convergence radii satisfy that $\Delta\geq\bar{\Delta}>0$. It indicates that after $s$ ($s\geq1$) epochs, the upper bound of $\mathbb{E}\left[\|\mathcal{G}_{\eta}(\widetilde{w}_s)\|^2\right]$ in (\ref{convergence2}) can be relatively lower than the one in (\ref{convergence1}). Therefore, it's theoretically supported that Algorithm \ref{alg2} has inspiring potential to outperform Algorithm \ref{alg1}. Essentially, the deterministic restart serves to remove any accumulated deviations occurred by stochastic optimization.

At start of each epoch, it's expected to compute the full gradient at $w_0^{(s)}$ in both Algorithm \ref{alg1} and Algorithm \ref{alg2} to derive the initial conjugate parameter $\beta_1$. In terms of complexity estimate, the deterministic restart does not affect computational overhead, but stores the computed full gradient and use it as the first search direction. Typically, the workload associated with storage is negligible.

Note that the condition (\ref{Thm1}) implies $\frac{2}{\eta_2}-L\gamma-3\geq0$ provided that $\eta_2\leq\frac{2}{3}$, then we have $\gamma\leq\frac{2-3\eta_2}{\eta_2L}$. If we define $\varpi(b)=\frac{(1+2\eta_2^2)(n-b)}{\eta_2b(n-1)}$, the suggestion is to choose $\gamma\leq \min\{\frac{2-3\eta_2}{\eta_2L}, \frac{-1+\sqrt{1-24m\eta_2\varpi(b)+16m\varpi(b)}}{4L\eta_2m\varpi(b)}\}$, which depends on $\sqrt m$. 

\paragraph{Further Discussion on the Gradient Dominant Scenario} Furthermore, nonsmooth composite functions may possess the gradient dominant property (Tran-Dinh et al., 2022 \cite{tran2022hybrid}), which can be defined as follows: 
\begin{definition}
	A composite function $F$ is $\tau_o$-gradient dominant, if it satisfies
	\begin{equation}
	\label{dominant}
	F(w)-F(w_{\star})\leq\frac{\tau_o}{2}\|{\mathcal{G}}_{\eta}(w)\|^{2},
	\end{equation}
for any $w\in\mathrm{dom}(F)$ and $\eta>0$, where $\tau_o>0$.	
\end{definition}

If the nonsmooth composite function $P$ in (\ref{problem}) is $\tau_o$-gradient dominant (\ref{dominant}), we are sufficient to obtain another linear convergence rate in the following Corollary \ref{corollary2}.

\begin{corollary}
\label{corollary2}
Assume that $P$ satisfies the gradient dominant condition (\ref{dominant}). Let Algorithm \ref{alg1} be implemented with $\widetilde{\eta}$ satisfying condition (\ref{Wolfe2}) and assume that during iterations $\widetilde{\eta}\in\{\eta_{1}, \eta_{2}\}$. Let
$$\delta'_k=\frac{\tau_o\xi}{2}\cdot(\mathbb{E}[\|\nabla f(\widetilde{w}_k)\|^2]+C), \quad k=0, 1, \cdots, s-1,$$
and $\delta'=\operatorname*{max}_{0\leq k\leq s-1}\delta'_{k}$. Assume that we choose $\xi'=\frac{\tau_o}{(m+1)\eta_1^2\gamma}<1$ and $b, \gamma$ that satisfy condition (\ref{Thm1}). For any tolerance $\varepsilon=\frac{\delta'}{1-\xi'}>0$, we have for Algorithm \ref{alg1} that
\begin{equation}
\mathbb{E}\left[P(\widetilde{w}_s)-P(w_{\star})\right]\leq(\xi')^s\left(\mathbb{E}\left[P(\widetilde{w}_0)-P(w_{\star})-\varepsilon\right]\right)+\varepsilon.
\end{equation}
Thus, $\left\{\mathbb{E}\left[P(\widetilde{w}_s)-P(w_{\star})\right]\right\}$ converges linearly with rate $\xi'$ to an $\varepsilon$-ball around zero.

\noindent
In Algorithm \ref{alg2} with $d_{0}^{(s)}=v_0^{(s)}$ ($s\geq1$), for any tolerance $\epsilon=\frac{\delta''}{1-\xi'}>0$, we have
\begin{equation}
\mathbb{E}\left[P(\widetilde{w}_s)-P(w_{\star})\right]\leq(\xi')^s\left(\mathbb{E}\left[P(\widetilde{w}_0)-P(w_{\star})-\epsilon\right]\right)+\epsilon,
\end{equation}
where we define $\delta''_k=\frac{\tau_o\xi}{2}\cdot(\mathbb{E}[\|\nabla f(\widetilde{w}_k)\|^2])$, $\delta''=\operatorname*{max}_{0\leq k\leq s-1}\delta''_{k}$.
Consequently, $\left\{\mathbb{E}\left[P(\widetilde{w}_s)-P(w_{\star})\right]\right\}$ converges linearly with rate $\xi'$ to an $\epsilon$-ball around zero.
\end{corollary}
\proof 
By Theorem \ref{theorem1}, we have
$$
\mathbb{E}\left[\|\mathcal{G}_{\eta}(w_{m}^{(s)})\|^{2}\right]
\leq\xi\left(\frac{(1-\hat{\beta})^2(1+\hat{\beta})}{\gamma\alpha(1+2\eta_2^2)\hat{\beta}^2}\mathbb{E}\left[P(w_{0}^{(s)})-P(w_{\star})\right]+\mathbb{E}\left[\|\nabla f(w_0^{(s)})\|^2\right]+C\right).
$$
By employing the gradient dominant condition (\ref{dominant}), one can show that
$$
\mathbb{E}\left[\|P(\widetilde{w}_s)-P(w_{\star})\right]
\leq\frac{\tau_o}{(m+1)\eta_1^2\gamma}\mathbb{E}\left[P(w_{0}^{(s)})-P(w_{\star})\right]+\frac{\tau_o \xi}{2}\left(\mathbb{E}\left[\|\nabla f(w_0^{(s)})\|^2\right]+C\right).
$$
By the definition of $\delta'$ and $\varepsilon=\frac{\delta'}{1-\xi'}$, it can be simplified as
$$
\mathbb{E}\left[P(\widetilde{w}_s)-P(w_{\star})\right]-\varepsilon\leq\xi'\left(\mathbb{E}\left[P(\widetilde{w}_{s-1})-P(w_{\star})-\varepsilon\right]\right).
$$
Therefore, we can show that
$$
\mathbb{E}\left[P(\widetilde{w}_s)-P(w_{\star})\right]\leq(\xi')^s\left(\mathbb{E}\left[P(\widetilde{w}_0)-P(w_{\star})-\varepsilon\right]\right)+\varepsilon.
$$
In Algorithm \ref{alg2}, we follow a similar line of reasoning and no longer expand it in detail.

\vspace{8pt}
It is worth noticing that the gradient dominant constant $\tau_o$ will determine the linear convergence constant $\xi'$ of $\left\{\mathbb{E}\left[P(\widetilde{w}_s)-P(w_{\star})\right]\right\}$.

\section{A Practical Variant}
We provide a practical variant in response to the concern that performing line searches (\ref{Wolfe1})-(\ref{Wolfe2}) can be time-consuming. To release computational (searching) workload, it's natural to consider reducing the number of calls to the inexact line search.

Indeed, it's not necessary to perform the line search in every step, especially in some later periods when iterates are well closed to the minimizer $w_{\star}$. Allocating additional computing power for (\ref{Wolfe1})-(\ref{Wolfe2}) may not be worthwhile all the time. Based on it, we present a class of practical variants that switch the search directions at regular or specific intervals and conduct the line search only in some essential steps. Recall from Section \ref{section} that the stochastic sufficient decrease condition (\ref{Wolfe1}) may not affect the convergence results, we omit the condition (\ref{Wolfe1}) to reduce the searching workload, then, the line search becomes 
\begin{equation}
\left|\langle v_{k+1}^{(s)}, d_{k+1-t}^{(s)}\rangle\right| \leq-c_{2}\langle v_{k+1-t}^{(s)}, d_{k+1-t}^{(s)}\rangle,  \label{Wolfe3}
\end{equation}
where $t$ ($t>1$) denotes the switching frequency. We perform (\ref{Wolfe3}) right before stochastic conjugate steps to ensure the mini-batch SARAH estimator $v_{k+1}^{(s)}$ (\ref{v}) used in the next stochastic conjugate step satisfies an essential property. For clarity, we show an example of regular intervals in Algorithm \ref{alg4}. Convergence analysis is presented after.

\subsection{A Practical Algorithm Using  Probabilistic Gradient}
We present a practical variant named as Acc-Prox-CG-SARAH-ST in Algorithm \ref{alg4}. Indeed, Algorithm \ref{alg4} utilizes a probabilistic gradient estimator (PAGE), which draws inspiration from a class of probabilistic switching approaches. Kovalev et al. (2020) \cite{kovalev2020switchSVRG} developed L-SVRG algorithm, which applies a probabilistic technique to switch between the GD estimator and the SVRG estimator during updates. In nonconvex optimization, Li et al. (2020a) \cite{li2020switchSARAH} chose to switch between GD and SARAH and proposed L2S method. Subsequently, Li et al. (2021) \cite{li2021switchpage} developed the PAGE estimator by switching mini-batch SGD with mini-batch SARAH, which allows flexible and adaptive probability scheme, as outlined in Algorithm \ref{alg3}.

\begin{algorithm}[h]
	\caption{ProbAbilistic Gradient Estimator (PAGE)}\label{alg3}
	\begin{algorithmic}[1]
		\Require initial point $w_0$, step size $\eta$, mini-batch size $b$ and $b'$ ($b'<b$), probability {$p_k\in(0, 1]$}
		\State $h^{0}=\frac{1}{b}\sum_{i\in \mathcal{B}}\nabla f_{i}(w^{0})$ // mini-batch size $\|\mathcal{B}\|=b$
		\For{$k = 0, 1,...$}
		\State $x^{k+1}=x^{k}-\eta h^{k}$
		\State $$h^{k+1}=\begin{cases}\frac{1}{b}\sum\limits_{i\in \mathcal{B}}\nabla f_i(w^{k+1})&\text{with probability} \quad p_k\\h^k+\frac{1}{b'}\sum\limits_{i\in \mathcal{B}'}(\nabla f_i(w^{k+1})-\nabla f_i(w^k))&\text{with probability}\quad 1-p_k\end{cases}$$
		\EndFor
		\State $\widehat{x}_T$ chosen uniformly from $\{x^t\}_{t\in[T]}$
	\end{algorithmic}
\end{algorithm}

In Algorithm \ref{alg4}, we switch stochastic conjugate gradients with mini-batch SARAH estimators (Nguyen et al., 2017 \cite{nguyen2017MBSARAH}) at a frequency of $t$ $(1<t\leq m-1)$. It's equivalent to choosing the search direction $h_{k}^s$ within the $s$-th epoch as follows:
$$
d_{k}^{(s)}=
\begin{cases}-v_{k}^{(s)}+\beta_kd_{k-t}^{(s)},  &\text{with probability} \quad p_k,\\
\nabla f_{\mathcal{B}_k}(w_k^{(s)})-\nabla f_{\mathcal{B}_k}(w_{k-1}^{(s)})+v_{k-1}^{(s)}, &\text{with probability}\quad 1-p_k,
\end{cases}
$$
where the probability indicator $\textbf{p}=(p_0, p_1, \cdots, p_{m-1})^T$ is a sparse binary vector with an $\ell_0$-norm $\|\textbf{p}\|_0=\lfloor \frac{m-1}{t} \rfloor$. Therefore, Algorithm \ref{alg4} shares similar properties with PAGE algorithm, L-SVRG method and L2S approach. Due to the switching frequency, we compute $\beta_k^{AFRs}$ by (will be explained later)
\begin{equation}
\label{AFR2}
\beta_{k}^{FRs} =\frac{\|v_{k}^{(s)}\|^{2}}{\|v_{k-t}^{(s)}\|^{2}}, \quad \beta_k^{AFRs}=\min \left\{\beta_o, \rho\beta_{k}^{FRs}\right\}.
\end{equation}
Similarly, we calculate $\beta_k^{FRPRs}$ by
\begin{equation}
\label{FRPR2}
\beta_{k}^{PRs}=\frac{\langle v_{k}^{(s)}, v_{k}^{(s)}-v_{k-t}^{(s)}\rangle}{\|v_{k-t}^{(s)}\|^2}, \quad \beta_k^{FRPRs}=\left\{\begin{array}{cc}
-\beta_k^{F Rs}, & \text { if } \beta_k^{P Rs}<-\beta_k^{F Rs}, \\
\beta_k^{P Rs}, & \text { if }\left|\beta_k^{P Rs}\right| \leq \beta_k^{F Rs}, \\
\beta_k^{F Rs} & \text { if } \beta_k^{P Rs}>\beta_k^{F Rs}.
\end{array}\right.
\end{equation}

\begin{algorithm}[h]
	\caption{Acc-Prox-CG-SARAH-ST}\label{alg4}
	\begin{algorithmic}[1]
		\Require initial point $\widetilde{w}_0$, update frequency $m$, mini-batch size $b$, momentum parameter $\gamma \in(0,1]$, threshold parameters $\eta_2>0$, $\beta_o>0$, $\rho>0$, switching frequency $1<t\leq m-1$, fixed step size $\eta_f>0$.
		\State $h^{(0)}=\nabla f(\widetilde{w}_0)$
		\For{$s = 1,2,...,S$}
		\State $w_0^{(s)}=\widetilde{w}_{s-1}$
		\State $v_0^{(s)}=\frac{1}{n} \sum_{i \in N} \nabla f_i(w_0^{(s)})=\nabla f(w_0^{(s)})$
		\State $d_0^{(s)}=-h^{(s-1)}$
		\State Set $\eta=\eta_f$
		\State $y_0^{(s)}=\operatorname{prox}_{\eta \varphi}(w_0^{(s)} + \eta d_0^{(s)})$
		\State $w_1^{(s)} =\left(1-\gamma\right) w_0^{(s)}+\gamma y_0^{(s)}$
		\For{$k = 1,2,...,m-1$}
		\State Pick mini-batch $\mathcal{B}_k \subset\{1, \ldots, n\}$ of size $b$ uniformly at random
		\State Update mini-batch SARAH estimator $v_{k}^{(s)}$ by
		\State $$v_{k}^{(s)}=\nabla f_{\mathcal{B}_k}(w_k^{(s)})-\nabla f_{\mathcal{B}_k}(w_{k-1}^{(s)})+v_{k-1}^{(s)}$$
		\If {$\mod(k, t)=0$}
		\State Calculate $\beta_k$ by AFRs (\ref{AFR2}) or FRPRs (\ref{FRPR2}) formula
		\State Compute the search direction by
		\State $d_k^{(s)}=-v_{k}^{(s)}+\beta_kd_{k-t}^{(s)}$
		\Else 
		\State Compute the search direction by
		\State $d_k^{(s)}=-v_{k}^{(s)}$
		\EndIf
		\If {$\mod(k+1, t)=0$}
		\State Call the line search and find $\widetilde{\eta}$ satisfying (\ref{Wolfe3})
		\State Set $\eta=\min\left\{\widetilde{\eta}, \eta_2 \right\}$
		\Else
		\State Set $\eta=\eta_f$
		\EndIf
		\State Update iterates by
		\State $y_{k}^{(s)}=\operatorname{prox}_{\eta\varphi}(w_k^{(s)} + \eta d_k^{(s)})$
		\State $w_{k+1}^{(s)} =\left(1-\gamma\right) w_k^{(s)}+\gamma y_{k}^{(s)}$
		\EndFor
		\State $h^{(s)}=v_m^{(s)}$
		\State $\widetilde{w}_s=w_m^{(s)}$
		\EndFor
	\end{algorithmic}
\end{algorithm}

Assume that $H$ denotes the set of all steps using stochastic conjugate directions (\ref{d}) with cardinality $q=|H|=\lfloor \frac{m-1}{t} \rfloor$, and $H'$ denotes the other complementary set containing other steps. Now, we display the related convergence results.
\begin{theorem}
	\label{theorem3}
    Let Algorithm \ref{alg4} be implemented with $\widetilde{\eta}$ satisfying condition (\ref{Wolfe3}) and switching frequency $1<t\leq m-1$. Assume that during iterations $\widetilde{\eta}\in\{\eta_{1}, \eta_{2}\}$. Parameters $b$, $\gamma$ are chosen under a suitable condition \footnote{This condition will be specified as (\ref{Thm3}) in Appendix D.}. Let $H$ denote the set of stochastic conjugate steps with $q=|H|=\lfloor \frac{m-1}{t} \rfloor$. For any $s>1$ in Algorithm \ref{alg4},
	$$
	\mathbb{E}\left[\|\mathcal{G}_{\eta}(w_{m}^{(s)})\|^{2}\right]
	\leq\frac{\xi^{st}(1-\hat{\beta})^2}{\gamma(1+2\eta_2^2)\hat{\beta}^2(1-\hat{\beta}^{q})^2}\mathbb{E}\left[P(w_{0}^{(s)})-P(w_{\star})\right]+\xi^{st}\mathbb{E}\left[\|\nabla f(w_0^{(s)})\|^2\right]+C^{st}\xi^{st},
	$$
	where parameters $\xi^{st}$ and $C^{st}$ are given by
	\begin{equation}
	\label{C2}
	\begin{aligned}
	&\xi^{st}=\frac{(2+4\eta_2^2)\alpha\hat{\beta}^2(1-\hat{\beta}^{q})}{(m+1)\eta_1^2(1-\hat{\beta})^2(1+\hat{\beta})},\\
	\quad
	&C^{st}=\left(\frac{(1-\hat{\beta})(1+\hat{\beta}^{q})}{\alpha}\tau+\frac{(1-\hat{\beta})^2(1+\hat{\beta})}{2\alpha\hat{\beta}^2(1-\hat{\beta}^{q})}\right)\sigma^2,
	\end{aligned}
	\end{equation}
\end{theorem}

\proof Deferred to the Appendix D.

\vspace{8pt}

\begin{theorem}
	\label{theorem4}
    Suppose that Theorem \ref{theorem3} holds. Let Algorithm \ref{alg4} be implemented with $\widetilde{\eta}$ satisfying condition (\ref{Wolfe3}) and switching frequency $1<t\leq m-1$. Assume that during iterations $\widetilde{\eta}\in\{\eta_{1}, \eta_{2}\}$. Let $H$ denote the set of stochastic conjugate steps with $q=|H|=\lfloor \frac{m-1}{t} \rfloor$. Define 
	$$\delta^{st}_k=\frac{(1-\hat{\beta})^2(1+\hat{\beta})}{\gamma\alpha(1+2\eta_2^2)\hat{\beta}^2(1-\hat{\beta}^{q})}\mathbb{E}\left[P(\widetilde{w}_k)-P(w_{\star})\right], \quad k=0, 1, \cdots, s-1,$$
	and $\delta^{st}=\operatorname*{max}_{0\leq k\leq s-1}\delta^{st}_{k}$. Assume that we choose $\xi^{st}=\frac{(2+4\eta_2^2)\alpha\hat{\beta}^2(1-\hat{\beta}^{q})}{(m+1)\eta_1^2(1-\hat{\beta})^2(1+\hat{\beta})}<1$ and $b$, $\gamma$ under condition (\ref{Thm3}). Then, by setting $\Delta^{st}=\frac{\xi^{st}(\delta^{st}+C)}{1-\xi^{st}}$, we have
	\begin{equation}
	\label{Thm4}
	\mathbb{E}\left[\|\mathcal{G}_{\eta}(\widetilde{w}_s)\|^2\right]-\Delta^{st}\leq(\xi^{st})^s\left(\|\nabla f(\widetilde{w}_0)\|^2-\Delta^{st}\right).
	\end{equation}
\end{theorem}
\proof  
It actually follows a similar line of reasoning as in the proof in Theorem \ref{theorem2}, here we will no longer elaborate it in detail.

\vspace{8pt}
The stochastic conjugate directions $\{d_k^{(s)}\}$ (\ref{d}) have demonstrated significant compatibility with $\{v_{k}^{(s)}\}$ (\ref{v}) in terms of convergence. The practical variant can also ensure stable and robust optimization, as indicated by (\ref{Thm4}), where the upper bound of $\mathbb{E}[\|\mathcal{G}_{\eta}(\widetilde{w}_s)\|^2]$ decreases as the epochs progress forward. The total number of calls to the line search (\ref{Wolfe3}) does not exceed $Sq$ times.

Actually, Algorithm \ref{alg4} shows an example of regular intervals. In addition, we can consider another scheme of switching the search directions only in later periods and at a low frequency (we call it the concentrated scheme). It entails regularly performing the line search in early periods to improve initial numerical performance. As $\{w_k^{(s)}\}$ approaches $w_{\star}$, we convert to predominantly utilizing $v_{k}^{(s)}$ with constant step sizes. Due to space limitations, we briefly introduce it and related analysis can be derived analogously as in Theorem \ref{theorem3} and Theorem \ref{theorem4}.

\section{Experiments}
In this section, we performed a series of experiments to verify the efficacy of our stochastic conjugate frameworks, including Algorithm \ref{alg1}, Algorithm \ref{alg2} and Algorithm \ref{alg4}. Particularly, we explore the properties of the algorithms by conducting different parametric settings. Then, we compare our algorithms with several state-of-art methods with fine-tuned parameters on four nonconvex sparse binary classification models.   All codes are available at \url{https://github.com/pzheng4218/SCG}.
\subsection{Experimental Settings }
We consider the (supervised) sparse binary classification model with a nonconvex loss function and a convex  $\ell_1$-regularizer as follows:
$$
\min_{w\in\mathbb{R}^d}\left\{P(w):=\frac{1}{n}\sum_{i=1}^{n}loss(a_{i}^{T}w,b_{i})+\lambda\|w\|_1\right\},
$$
where $\{(a_i,b_i)\}_{i=1}^n\subset\mathbb{R}^d\times\{-1, 1\}^n$ denotes the sequence of $n$ labeled pairs. Here, $\varphi=\lambda\|w\|_1$ is commonly referred to as the soft-thresholding (Donoho et al., 1994 \cite{donoho1994softthreholding}). In the experiments, we evaluate algorithms based on four nonconvex loss functions. The first one is chosen from Metel et al. (2019) \cite{metel2019loss1}, while the other three have been discussed in Zhao et al., (2010) \cite{zhao2010loss2}. Detailed descriptions are provided as follows:

\vspace{4pt}
\noindent
\textbf{(1) Lorenz loss $\ell_a$} 
\begin{equation}
\label{loss1}
	loss(a_{i}^{T}w,b_{i})=\ell_a:=\left\{\begin{array}{cc}
	\ln(1+(b_ia_i^Tw-1)^2) & \quad \text { if } b_ia_i^Tw\leq 1 \\
	0 & \quad \text { otherwise }
	\end{array}\right.
\end{equation}
We have the gradient with respect to $w$ as:
$$
\nabla_w \ell_a:=\left\{\begin{array}{cc}
\frac{2(b_ia_i^Tw-1)b_i}{1+(b_ia_i^Tw-1)^2}a_i & \quad \text { if } b_ia_i^Tw\leq 1 \\
0 & \quad \text { otherwise }
\end{array}\right.
$$
Lorenz loss function $\ell_a$ is $L$-smooth with Lipschitz constant $L = 4$.

\vspace{4pt}
\noindent
\textbf{(2) Normalized sigmoid loss $\ell_b$}
\begin{equation}
\label{loss2}
loss(a_{i}^{T}w,b_{i})=\ell_b:= 1-tanh(b_ia_{i}^{T}w)=\frac{2\exp\{-b_ia_i^Tw\}}{\exp\{b_ia_i^Tw\}+\exp\{-b_ia_i^Tw\}}
\end{equation}
The gradient of (\ref{loss2}) with respect to $w$ is:
$$
\nabla_w \ell_b:= \frac{-4b_i}{\left(\exp\{b_ia_i^Tw\}+\exp\{-b_ia_i^Tw\}\right)^2}a_i
$$
Normalized sigmoid loss function $\ell_b$ is $L$-smooth with Lipschitz constant $L\approx0.7698$.

\vspace{4pt}
\noindent
\textbf{(3) Logistic difference loss $\ell_c$}
\begin{equation}
\label{loss3}
loss(a_{i}^{T}w,b_{i})=\ell_c:= \ln(1+\exp\{-b_ia_i^Tw\})-\ln(1+\exp\{-b_ia_i^Tw-1\})
\end{equation}
Then, we have its gradient with respect to $w$ as:
$$
\nabla_w \ell_c:= \frac{-4b_i}{\left(\exp\{b_ia_i^Tw\}+\exp\{-b_ia_i^Tw\}\right)^2}a_i
$$
This function $\ell_c$ is also $L$-smooth with $L=0.092372$.

\vspace{4pt}
\noindent
\textbf{(4) Nonconvex loss in $2$-layer neural network $\ell_d$}
\begin{equation}
\label{loss4}
loss(a_{i}^{T}w,b_{i})=\ell_d:= \left(1 - \frac{1}{1+\exp\{-b_ia_i^Tw\}}\right)^2
\end{equation}
We have the gradient of (\ref{loss4}) with respect to $w$ as:
$$
\nabla_w \ell_d:=\frac{-2b_i\exp^2\{-b_ia_i^Tw\}}{(1+\exp\{-b_ia_i^Tw\})^3}a_i
$$
Note that the gradient $\nabla_w \ell_d$ is $L$-Lipshitz continuous with $L=0.15405$.

We evaluate algorithms across distinct $\ell_a$, $\ell_b$, $\ell_c$ and $\ell_d$ to verify their scalability and robustness in practical implementations.

Datasets used in the experiments are publicly available at LIBSVM \footnote{ \url{https://www.csie.ntu.edu.tw/~cjlin/libsvmtools/datasets/}.}online. Specifically, we utilize three representative datasets: $w8a$ ($n= 49,749, d=300$), $a9a$ ($n= 32,561, d=123$), and $gisette$ ($n= 6,000, d=5,000$). During preprocessing operations, we scale the data matrix $X$ for all dimensions into range of $[-1, 1]$ by $\ell_2$-normalization. The $\ell_1$-regularization coefficient $\lambda$ is set to be $\lambda=10^{-2}/n\approx2\times10^{-7}$ for $w8a$ data set; $\lambda=10^{-3}/n\approx3\times10^{-8}$ for $a9a$ data set; and $\lambda=10^{-7}\sqrt{n}/\sqrt{d}\approx1.1\times10^{-7}$ for $gisette$ data set.

To ensure fair comparisons, we compute the $\ell_2$-norm of the gradient mapping $\{\|\mathcal{G}_\eta(w_k^{(s)})\|^2\}$ under the same value of $\eta=0.5$ for all algorithms. To compute the sub-optimality $\{P(\widetilde{w}_{s})-P(w_{\star})\}$, we derive the optimal loss $P(w_{\star})$ through all algorithms.

\subsection{Properties of Algorithms}
As observed from Algorithm \ref{alg1} and Algorithm \ref{alg2}, the numerical efficiency of performance primarily rely on the mini-batch size $b$, the epoch length $m$, the conjugate formula $\beta$ and the momentum parameter $\gamma$. In this subsection, we randomly select the normalized sigmoid loss $\ell_b$ (\ref{loss2}) model to explore basic properties of our algorithms. For illustration, we display the parametric configurations of Algorithm \ref{alg1} and Algorithm \ref{alg2} in the following Table \ref{table2}.
\begin{table}
	\setlength{\abovecaptionskip}{0pt}
	\setlength{\belowcaptionskip}{1pt}
	\centering
	\vspace{4pt}
	\renewcommand\arraystretch{1.4}
	\resizebox{0.9\textwidth}{!}{
	\begin{tabular}{lccccc}
		\hline
		\multicolumn{1}{l}{Algorithms}&\multicolumn{1}{l}{\quad $b$} &\multicolumn{1}{l}{\quad $m$} &\multicolumn{1}{c}{$\beta$} &\multicolumn{1}{c}{$\gamma$} &\multicolumn{1}{c}{abbrev.}
		\\ \hline\noalign{\smallskip}
		Acc-Prox-CG-SARAH (-RS) & $\lfloor n^{\frac{1}{3}} \rfloor$  & $\lfloor \frac{1}{3}n^{\frac{1}{3}} \rfloor$ & $AFR$ (\ref{AFR}) & $\frac{\sqrt{m}}{4}$ & $v_1$ ($RS-v_1$)\\
		Acc-Prox-CG-SARAH (-RS) & $\lfloor n^{\frac{1}{3}} \rfloor$  & $\lfloor \frac{1}{3}n^{\frac{1}{3}} \rfloor$ & $FRPR$ (\ref{FRPR}) & $\frac{\sqrt{m}}{4}$ & $v_2$ ($RS-v_2$)\\
		Acc-Prox-CG-SARAH (-RS)& $\lfloor \sqrt{n} \rfloor$  & $\lfloor \frac{1}{3}n^{\frac{1}{3}} \rfloor$ & $AFR$ (\ref{AFR}) & $\frac{\sqrt{m}}{4}$ & $v_3$ ($RS-v_3$)\\
		Acc-Prox-CG-SARAH (-RS)& $\lfloor \sqrt{n} \rfloor$  & $\lfloor \frac{1}{3}n^{\frac{1}{3}} \rfloor$ & $AFR$ (\ref{AFR}) & $\frac{\sqrt{m}}{5}$ & $v_4$ ($RS-v_4$)\\
		Acc-Prox-CG-SARAH (-RS)& $\lfloor \sqrt{n} \rfloor$  & $\lfloor \frac{1}{3}n^{\frac{1}{3}} \rfloor$ & $FRPR$ (\ref{FRPR}) & $\frac{\sqrt{m}}{4}$ & $v_5$ ($RS-v_5$)\\
		Acc-Prox-CG-SARAH (-RS)& $\lfloor \sqrt{n} \rfloor$  & $\lfloor \frac{1}{3}n^{\frac{1}{3}} \rfloor$ & $FRPR$ (\ref{FRPR}) & $\frac{\sqrt{m}}{5}$ & $v_6$ ($RS-v_6$)\\
		Acc-Prox-CG-SARAH (-RS)& $\lfloor 2\sqrt{n} \rfloor$  & $\lfloor \frac{1}{3}n^{\frac{1}{3}} \rfloor$ & $AFR$ (\ref{AFR}) & $\frac{\sqrt{m}}{4}$ & $v_7$ ($RS-v_7$)\\
		Acc-Prox-CG-SARAH (-RS)& $\lfloor 2\sqrt{n} \rfloor$  & $\lfloor \frac{1}{3}n^{\frac{1}{3}} \rfloor$ & $FRPR$ (\ref{FRPR}) & $\frac{\sqrt{m}}{4}$ & $v_8$ ($RS-v_8$)\\
		\hline
	\end{tabular}}
    \caption{A description of parameter configurations: for brevity, we use abbrev. $v_1$ to $v_8$ to denote different settings of Algorithm \ref{alg1}, and abbrev. $RS-v_1$ to $RS-v8$ to denote the corresponding settings of Algorithm \ref{alg2}. For fair comparisons, we unify the epoch length $m:=\lfloor n^{\frac{1}{3}} \rfloor$. According to Theorem \ref{theorem3}, we configure the momentum parameter $\gamma$ relying on $\sqrt{m}$.} 
	\label{table2}
\end{table}

\begin{figure*}[htbp]
	\centering
	\subfigure[w8a]
	{\includegraphics[width=0.327\textwidth]{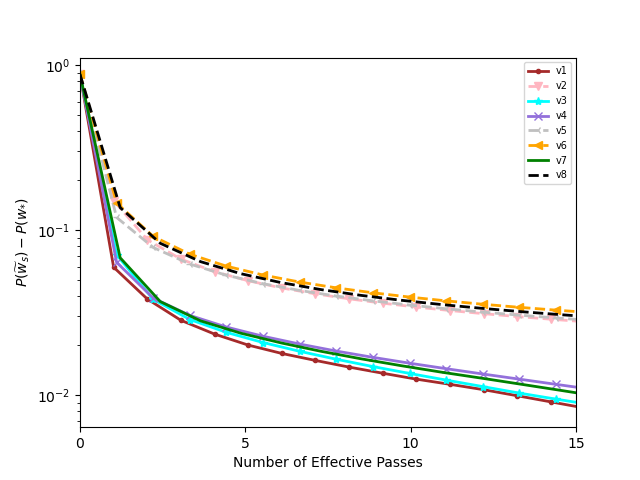}}
	\subfigure[a9a]
	{\includegraphics[width=0.327\textwidth]{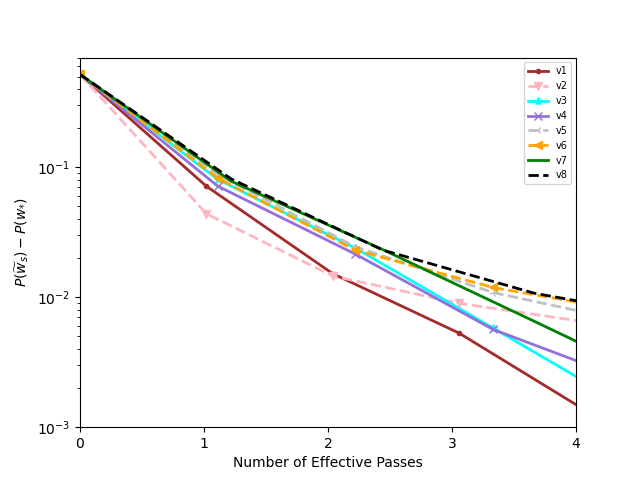}}
	\subfigure[gisette]
	{\includegraphics[width=0.327\textwidth]{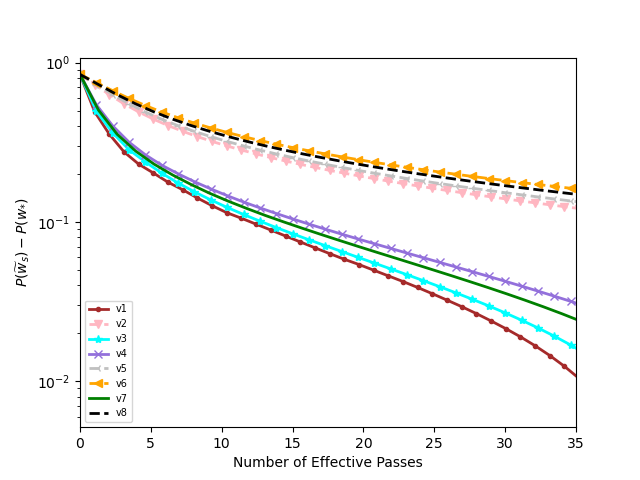}}
	
	\subfigure[w8a]
	{\includegraphics[width=0.327\textwidth]{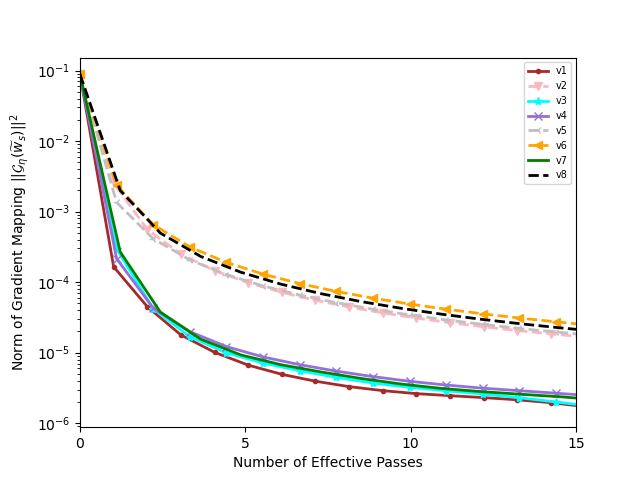}}
	\subfigure[a9a]
	{\includegraphics[width=0.327\textwidth]{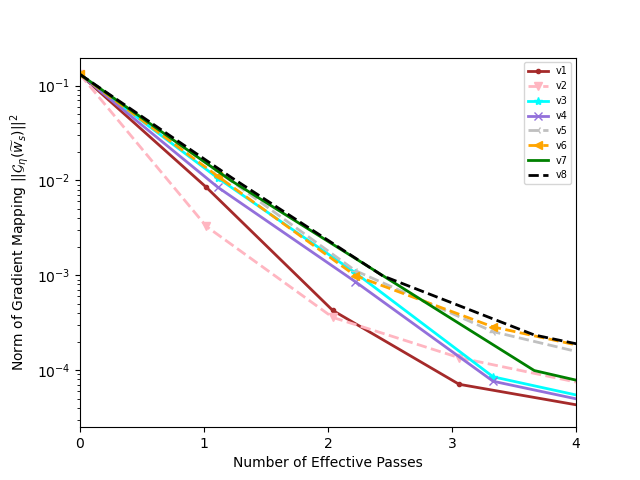}}
	\subfigure[gisette]
	{\includegraphics[width=0.327\textwidth]{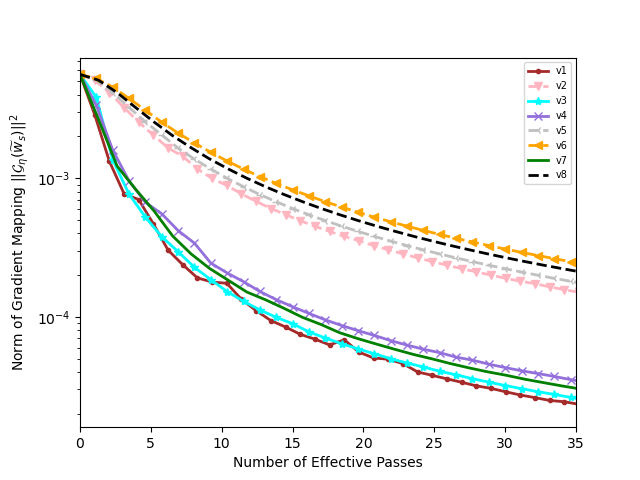}}
	\caption{\footnotesize Comparisons of different variants of Algorithm \ref{alg1}, including $v_1$, $v_2$, $v_3$, $v_4$, $v_5$, $v_6$, $v_7$ and $v_8$.}
	\label{fig1}
\end{figure*}

First, we implement $8$ variants of Algorithm \ref{alg1} ($v_1$ to $v_8$) on data sets $w8a$, $a9a$ and $gisette$, respectively. The associated results are shown in Fig. \ref{fig1}, where the $x$-axis denotes the number of effective passes and the $y$-axis stands for the sub-optimality $\{P(\widetilde{w}_{s})-P(w_{\star})\}$ or  $\{||\mathcal{G}_{\eta}(\widetilde{w}_{s})||^2\}$. From the results of $\{v_1, v_2, v_3, v_4, v_5, v_6\}$: we view that the variants in AFR (\ref{AFR}) updates ($v_1, v_3, v_4$) perform much faster than the others in FRPR (\ref{FRPR}) updates ($v_2, v_5, v_6$) in general. Still, it is worth noting that $v_2$ outperforms $v_1$ in certain early stages in $a9a$, but its sub-optimality stagnates around $10^{-2}$ accuracy during later optimization. Hence, to achieve the best performance, a pre-search over the conjugate updates is recommended for each particular data set. From the results of $\{v_1, v_2, v_3, v_5, v_7, v_8\}$: we observe that the variants with smaller mini-batch sizes ($v_1, v_2$) consistently work well. Thus, it's not necessary to apply a large batch in our Algorithm \ref{alg1}, as it may also require additional computing power. From the results of $\{v_3, v_4, v_5, v_6\}$: it suggests that choosing $\gamma$ closer to $1$ is a more preferable and advantageous choice to accelerate our Algorithm \ref{alg1}. 

\begin{figure*}[htbp]
	\centering
	\subfigure[w8a]
	{\includegraphics[width=0.327\textwidth]{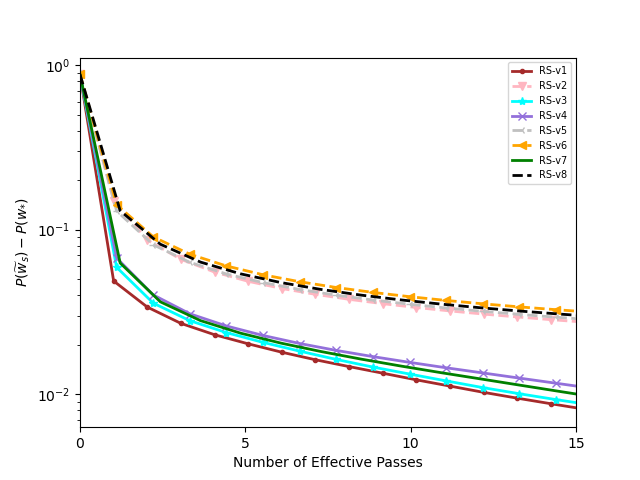}}
	\subfigure[a9a]
	{\includegraphics[width=0.327\textwidth]{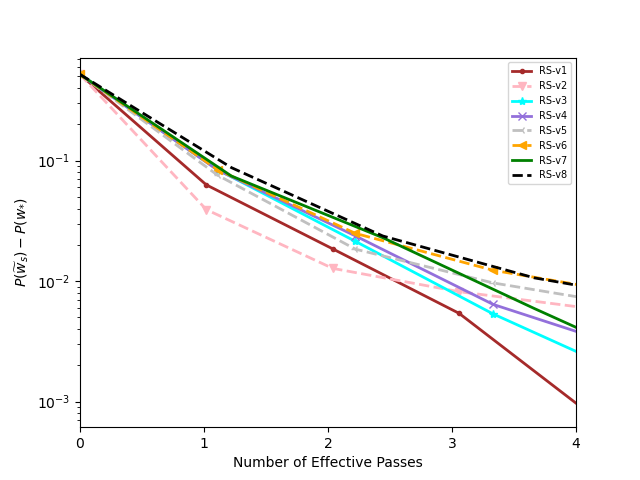}}
	\subfigure[gisette]
	{\includegraphics[width=0.327\textwidth]{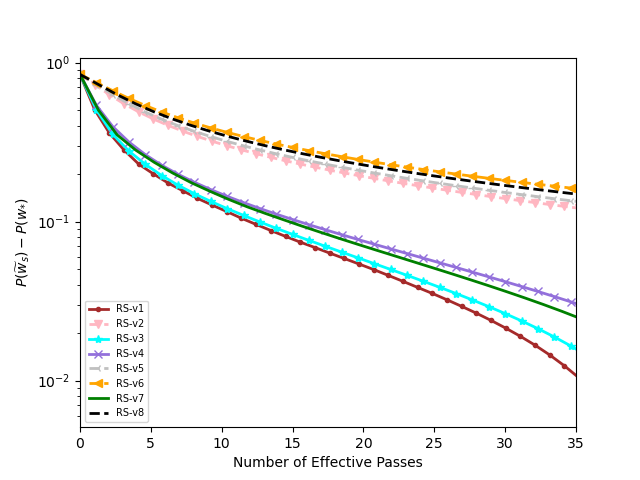}}
	
	\subfigure[w8a]
	{\includegraphics[width=0.327\textwidth]{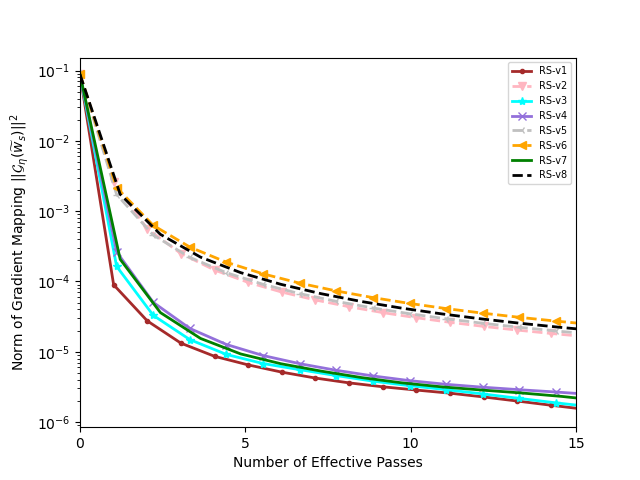}}
	\subfigure[a9a]
	{\includegraphics[width=0.327\textwidth]{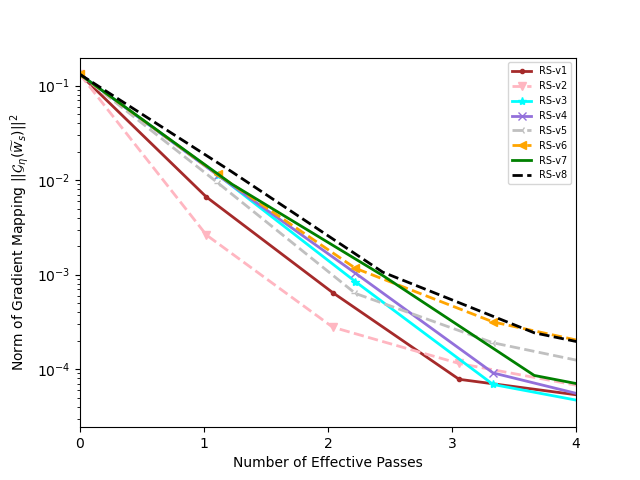}}
	\subfigure[gisette]
	{\includegraphics[width=0.327\textwidth]{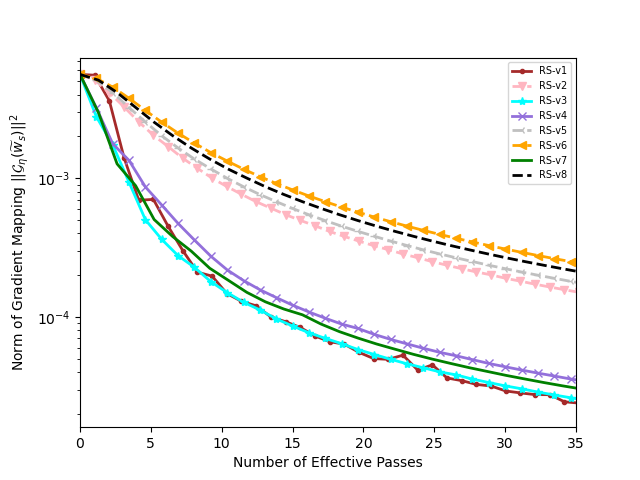}}
	\caption{\footnotesize Comparisons of different variants of Algorithm \ref{alg2}, including $RS-v_1$, $RS-v_2$, $RS-v_3$, $RS-v_4$, $RS-v_5$, $RS-v_6$, $RS-v_7$ and $RS-v_8$.}
	\label{fig2}
\end{figure*}

Second, we perform parallel tests on the variants of Algorithm \ref{alg2} ($RS-v_1$ to $RS-v_8$), as illustrated in Fig. \ref{fig2}. The results exhibit a substantial resemblance to those in Fig. \ref{fig1}. These two algorithms share similar properties to some great extend.

\begin{figure*}[htbp]
	\centering
	\subfigure[w8a]
	{\includegraphics[width=0.327\textwidth]{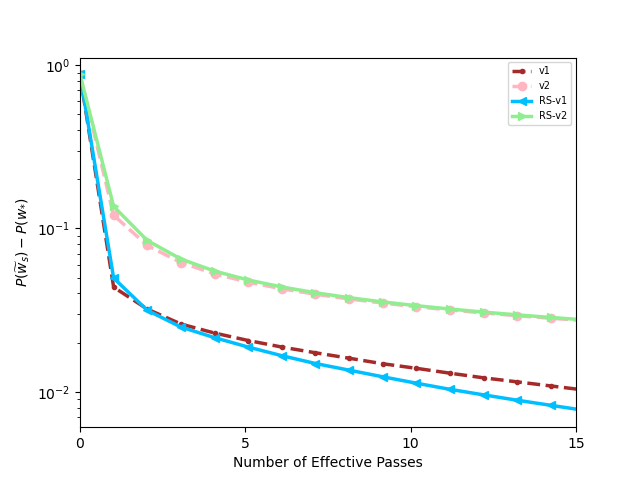}}
	\subfigure[a9a]
	{\includegraphics[width=0.327\textwidth]{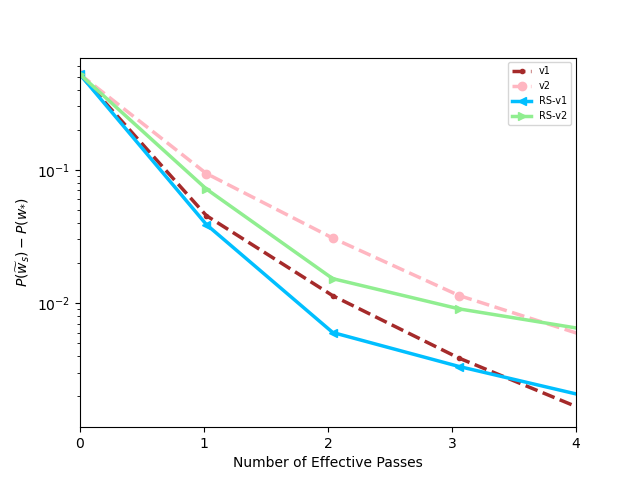}}
	\subfigure[gisette]
	{\includegraphics[width=0.327\textwidth]{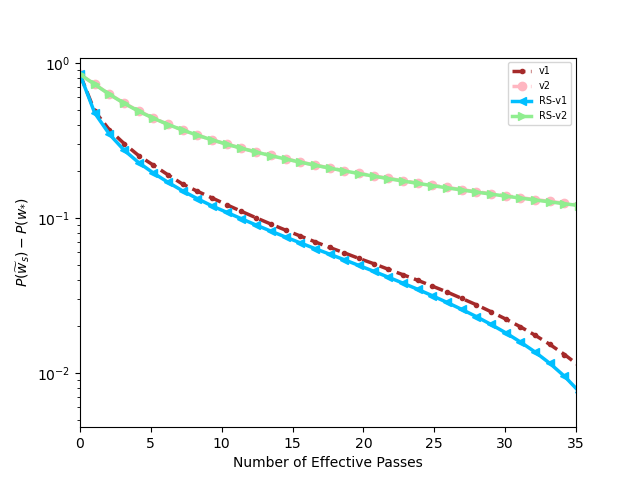}}
	
	\subfigure[w8a]
	{\includegraphics[width=0.327\textwidth]{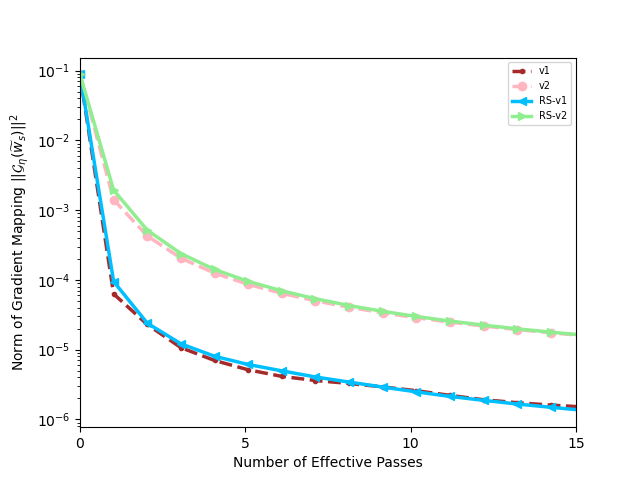}}
	\subfigure[a9a]
	{\includegraphics[width=0.327\textwidth]{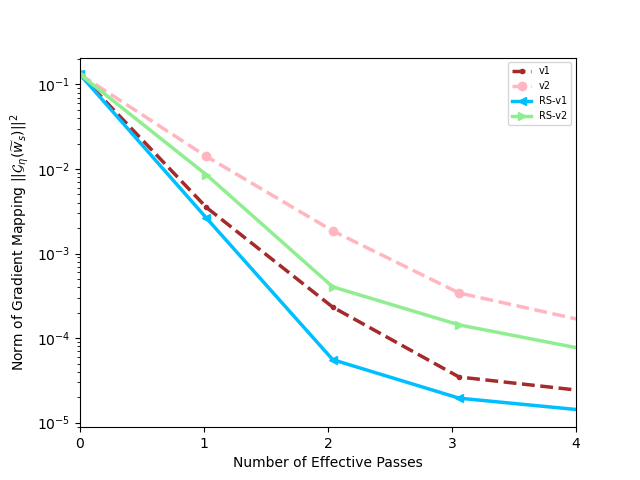}}
	\subfigure[gisette]
	{\includegraphics[width=0.327\textwidth]{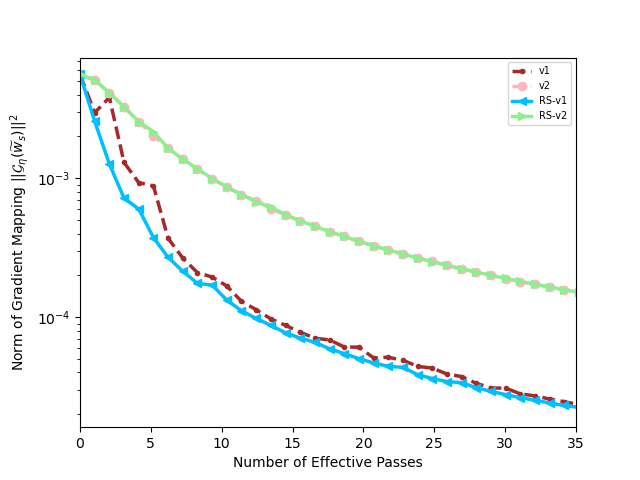}}
	\caption{\footnotesize Comparisons between $v_1$ and $RS-v_1$, $v_2$ and $RS-v_2$.}
	\label{fig3}
\end{figure*}

Eventually, we select the variants $v_1$ and $RS-v_1$ that achieve the best performance in AFR (\ref{AFR}) updates, along with $v_2$ and $RS-v_2$ which perform the best in FRPR (\ref{FRPR}) updates, to validate the efficacy of our deterministic restart scheme. The results are plotted in Fig. \ref{fig3}. As can be observed, under AFR (\ref{AFR}) updates, $RS-v_1$ consistently perform better across all the data sets, which is consistent with our theoretical analysis before. However, under FRPR (\ref{FRPR}) updates, the performance of $RS-v_2$ is comparable with $v_2$ in both $w8a$ and $gisette$. Overall, the deterministic restart scheme is more sensitive to algorithms using AFR (\ref{AFR}) update formula.

\subsection{Comparisons with other state-of-art methods}
In this subsection, we compare the convergence of Algorithm \ref{alg1} and Algorithm \ref{alg2} with other state-of-art methods across all the nonconvex and nonsmooth model, including $\ell_a$, $\ell_b$. $\ell_c$ and $\ell_d$. Specifically, we apply the configurations of $v_1$ and $RS-v_1$ while fine-tune each of the other algorithms according to their original literature, which are summarized in the Table \ref{table3}.

\begin{table}[h]
	\setlength{\abovecaptionskip}{0pt}
	\setlength{\belowcaptionskip}{1pt}
	\centering
	\vspace{4pt}
	\renewcommand\arraystretch{1.5}
	\resizebox{1.1\textwidth}{!}{
		\begin{tabular}{lcccc}
			\hline
			\multicolumn{1}{l}{Algorithms}&\multicolumn{1}{c}{mini-batch sizes and their related parameters} &\multicolumn{1}{c}{$m$} &\multicolumn{1}{c}{$\eta$} &\multicolumn{1}{c}{$\gamma$} 
			\\ \hline\noalign{\smallskip}
			Acc-Prox-CG-SARAH & $b=\lfloor n^{\frac{1}{3}} \rfloor$  & $\lfloor \frac{1}{3}n^{\frac{1}{3}} \rfloor$ & (\ref{Wolfe2}) & $\frac{\sqrt{m}}{4}$ \\
			Acc-Prox-CG-SARAH-RS & $b=\lfloor n^{\frac{1}{3}} \rfloor$  & $\lfloor \frac{1}{3}n^{\frac{1}{3}} \rfloor$ & (\ref{Wolfe2}) & $\frac{\sqrt{m}}{4}$\\
			ProxSARAH \cite{phamproxsarah} & $C=\frac{2}{3L^{2}\gamma^{2}}, b=\lfloor\frac{n^{\frac{1}{3}}}{C}\rfloor$  & $\lfloor n^{\frac{1}{3}}\rfloor$ & $\frac{2}{4+L\gamma}$ & $0.99$\\
			(Prox-)Spiderboost \cite{wang2019spiderboost} & $b=\lfloor \sqrt{n} \rfloor$  & $\lfloor \sqrt{n} \rfloor$ & $\frac{1}{2L}$ & / \\
			ProxSVRG+ \cite{li2018proxsvrg+} & $B=\lfloor \frac{n}{5} \rfloor, b=\lfloor n^{\frac{2}{3}} \rfloor$  & $\lfloor \sqrt{b} \rfloor$ & $\frac{1}{6L}$ & /\\
			ProxHSGD-RS \cite{tran2022hybrid} & $\hat{b}=b, \tilde{b}=c_1^2[b(m+1)]^{1/3}, c_{1}=\frac{b^{1/3}}{(m+1)^{2/3}}, \beta=1-\frac{\hat{b}^{1/2}}{[\tilde{b}(m+1)]^{1/2}}$  & $\lfloor n^{\frac{1}{3}} \rfloor$ & $\frac{1}{L}$ & $0.95$\\
			\hline
	\end{tabular}}
	\caption{A description of parameter configurations of all methods: ProxSARAH is configured following section 3.4 in \cite{phamproxsarah}; (Prox-)Spiderboost is established according to Theorem 1 in \cite{wang2019spiderboost}; ProxSVRG+ is implemented in accordance with Theorem 1 in \cite{li2018proxsvrg+}; ProxHSGD-RS is configured according to Theorem 3 in \cite{tran2022hybrid}, where we choose the upper bound $\eta =\frac{1}{L}$ for full efficiency.} 
	\label{table3}
\end{table}

\begin{figure*}[h]
	\centering
	\subfigure[w8a]
	{\includegraphics[width=0.327\textwidth]{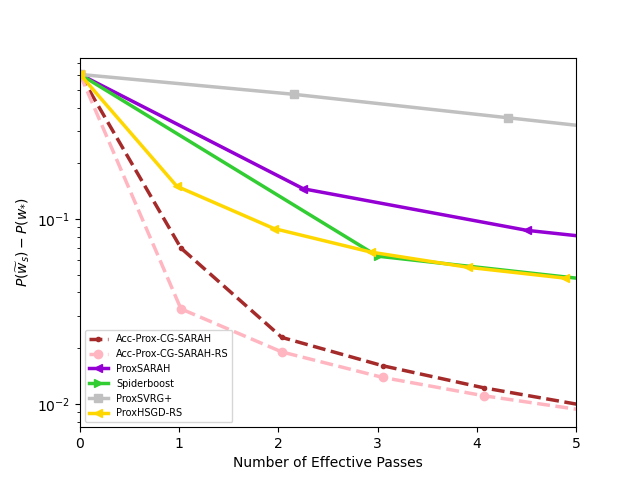}}
	\subfigure[a9a]
	{\includegraphics[width=0.327\textwidth]{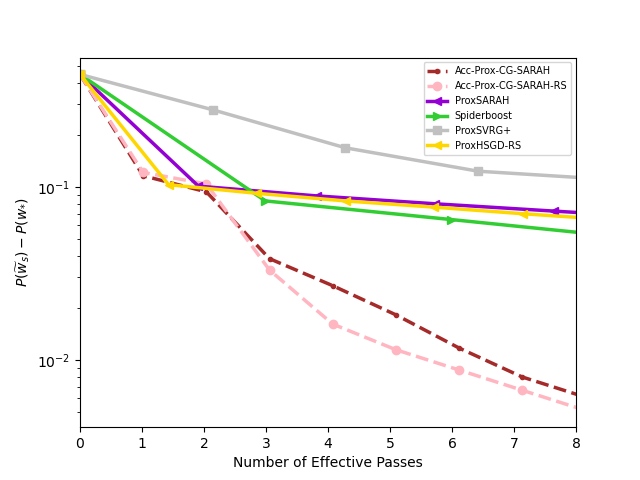}}
	\subfigure[gisette]
	{\includegraphics[width=0.327\textwidth]{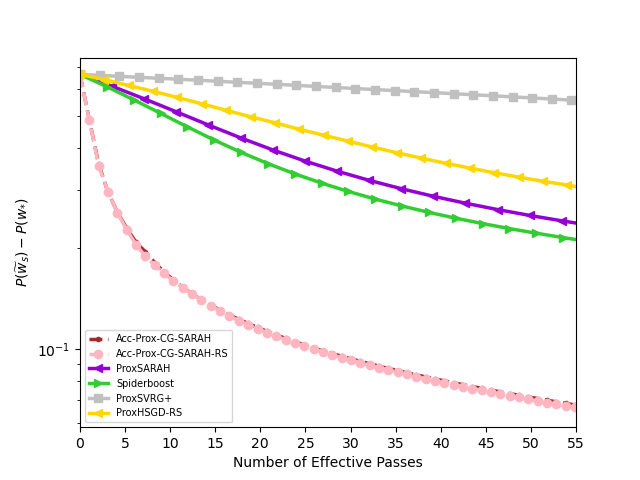}}
	
	\subfigure[w8a]
	{\includegraphics[width=0.327\textwidth]{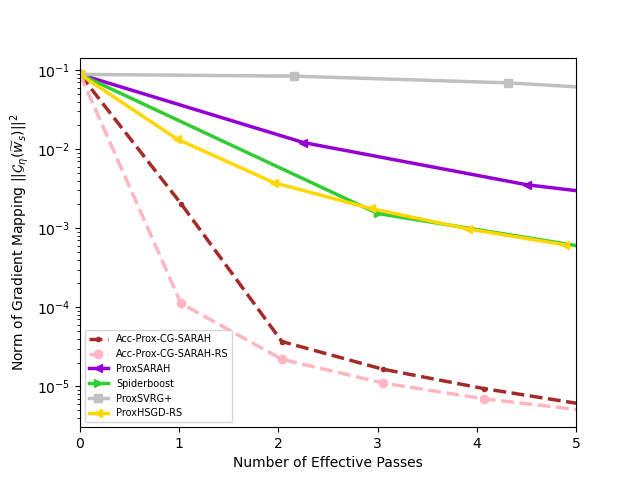}}
	\subfigure[a9a]
	{\includegraphics[width=0.327\textwidth]{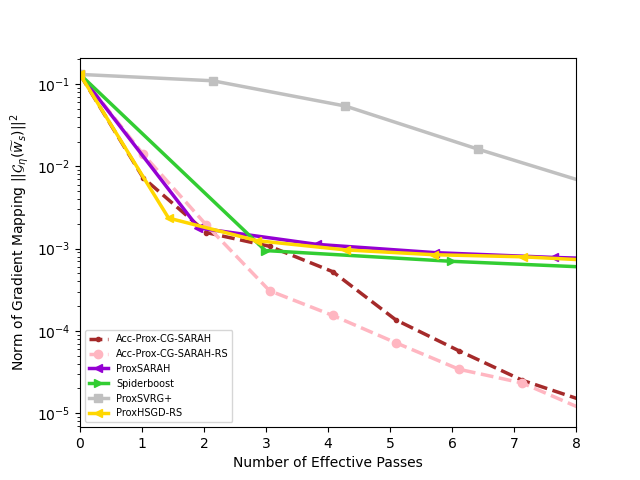}}
	\subfigure[gisette]
	{\includegraphics[width=0.327\textwidth]{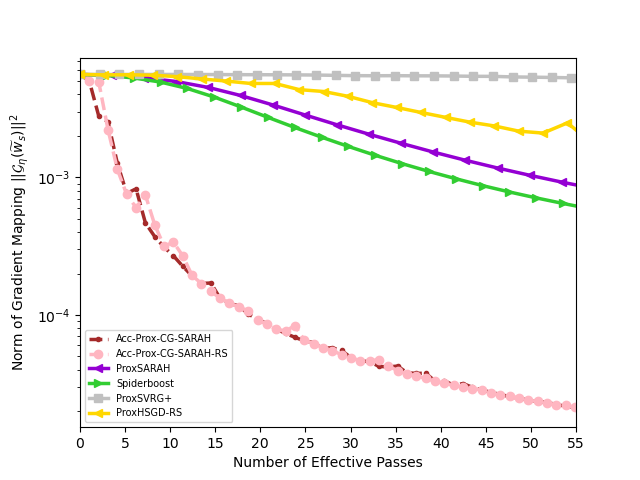}}
	\caption{\footnotesize Comparisons with other state-of-art methods over Lorenz loss $\ell_a$.}
	\label{fig4}
\end{figure*}

\begin{figure*}[h]
	\centering
	\subfigure[w8a]
	{\includegraphics[width=0.327\textwidth]{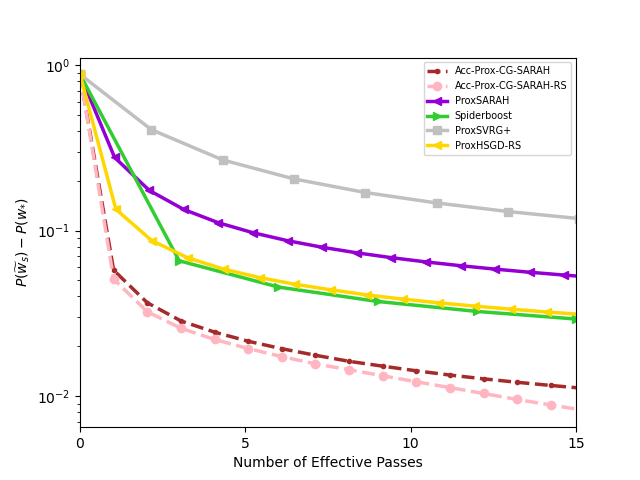}}
	\subfigure[a9a]
	{\includegraphics[width=0.327\textwidth]{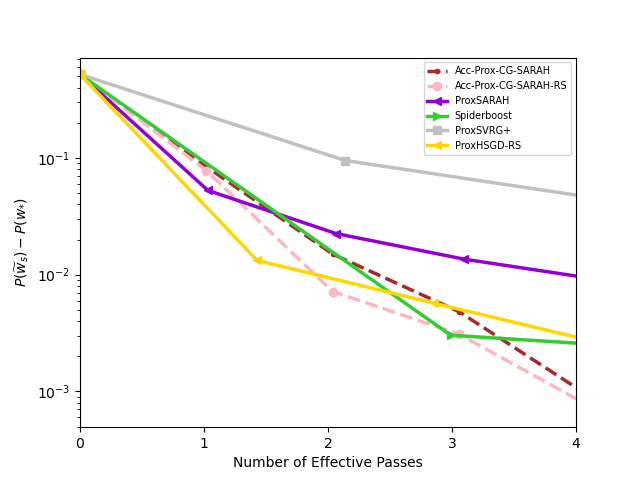}}
	\subfigure[gisette]
	{\includegraphics[width=0.327\textwidth]{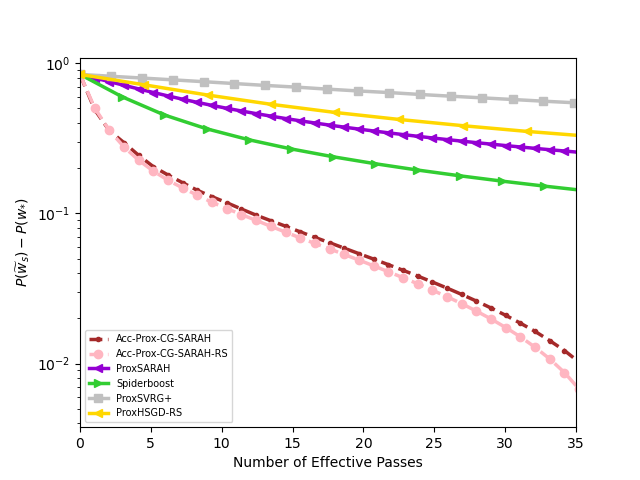}}
	
	\subfigure[w8a]
	{\includegraphics[width=0.327\textwidth]{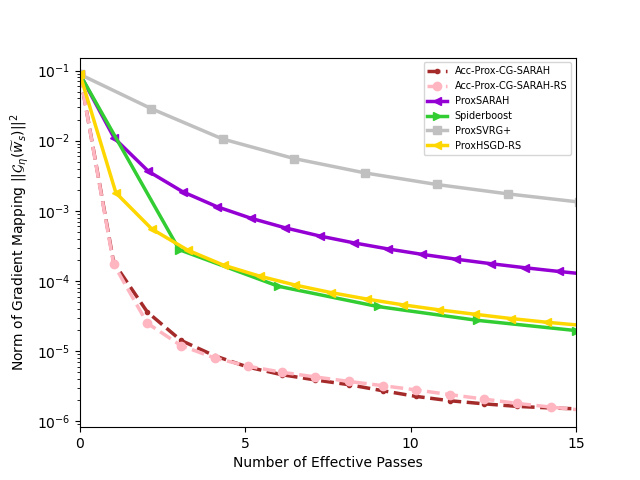}}
	\subfigure[a9a]
	{\includegraphics[width=0.327\textwidth]{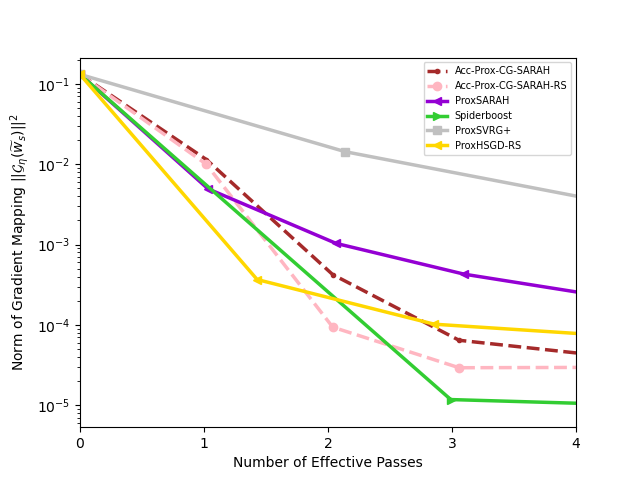}}
	\subfigure[gisette]
	{\includegraphics[width=0.327\textwidth]{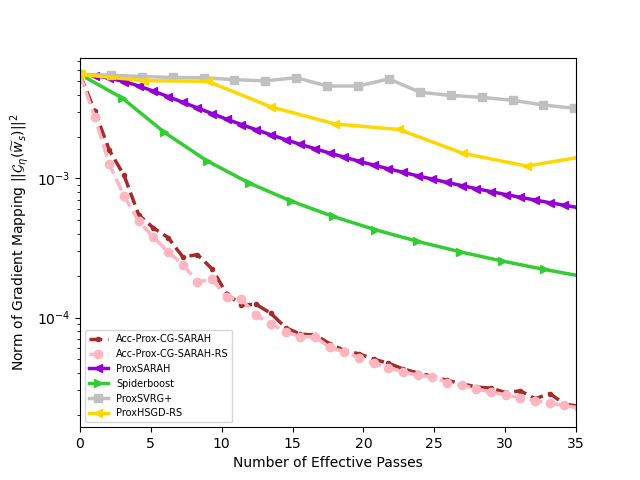}}
	\caption{\footnotesize Comparisons with other state-of-art methods over normalized sigmoid loss $\ell_b$.}
	\label{fig5}
\end{figure*}

\begin{figure*}[h]
	\centering
	\subfigure[w8a]
	{\includegraphics[width=0.327\textwidth]{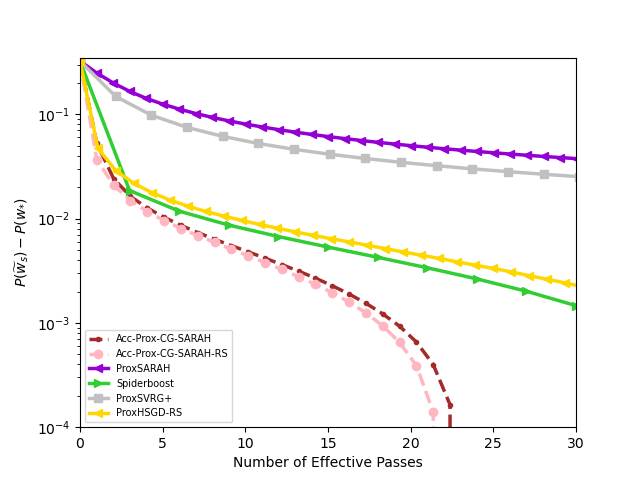}}
	\subfigure[a9a]
	{\includegraphics[width=0.327\textwidth]{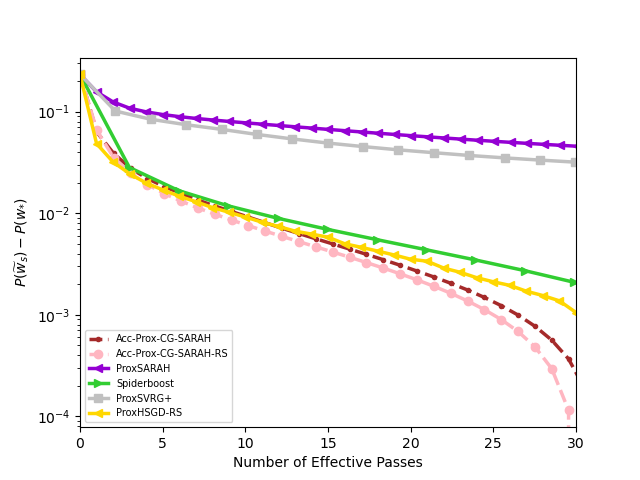}}
	\subfigure[gisette]
	{\includegraphics[width=0.327\textwidth]{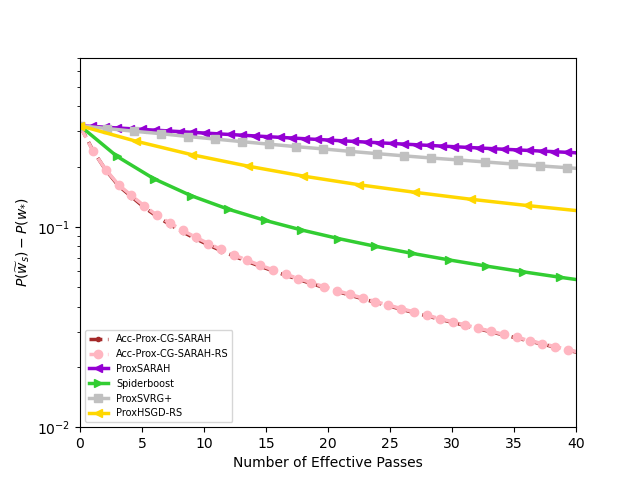}}
	
	\subfigure[w8a]
	{\includegraphics[width=0.327\textwidth]{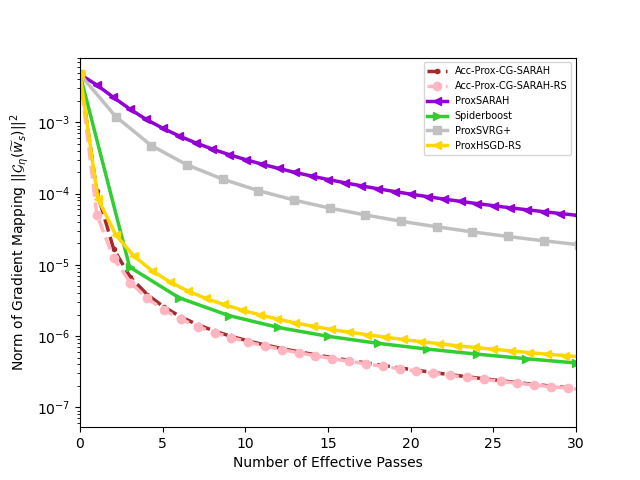}}
	\subfigure[a9a]
	{\includegraphics[width=0.327\textwidth]{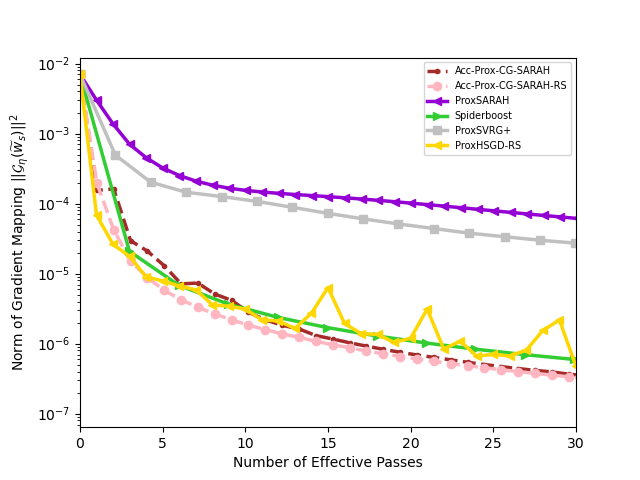}}
	\subfigure[gisette]
	{\includegraphics[width=0.327\textwidth]{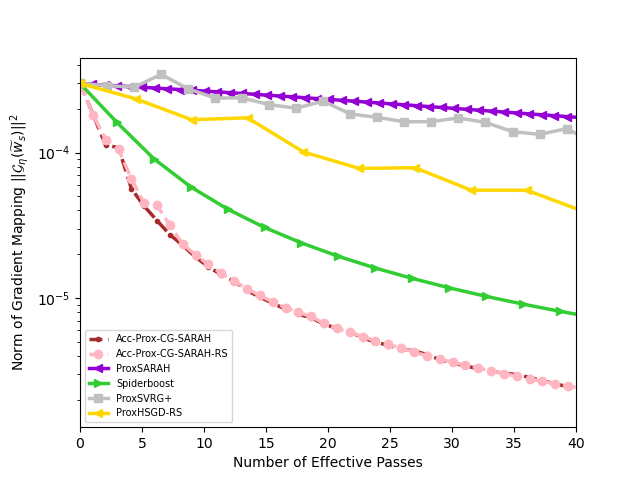}}
	\caption{\footnotesize Comparisons with other state-of-art methods over logistic difference loss $\ell_c$.}
	\label{fig6}
\end{figure*}

\begin{figure*}[h]
	\centering
	\subfigure[w8a]
	{\includegraphics[width=0.327\textwidth]{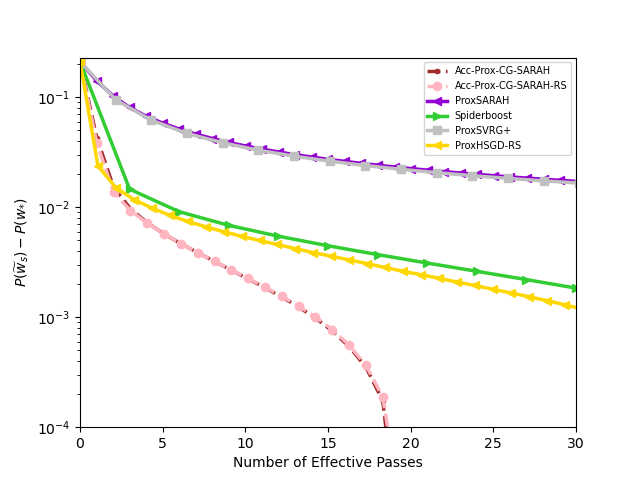}}
	\subfigure[a9a]
	{\includegraphics[width=0.327\textwidth]{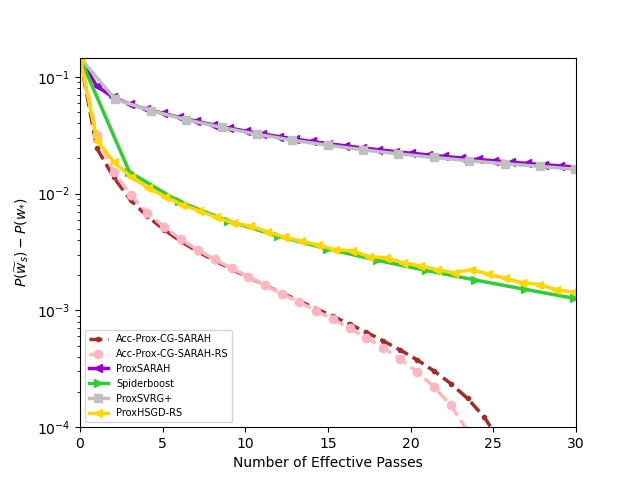}}
	\subfigure[gisette]
	{\includegraphics[width=0.327\textwidth]{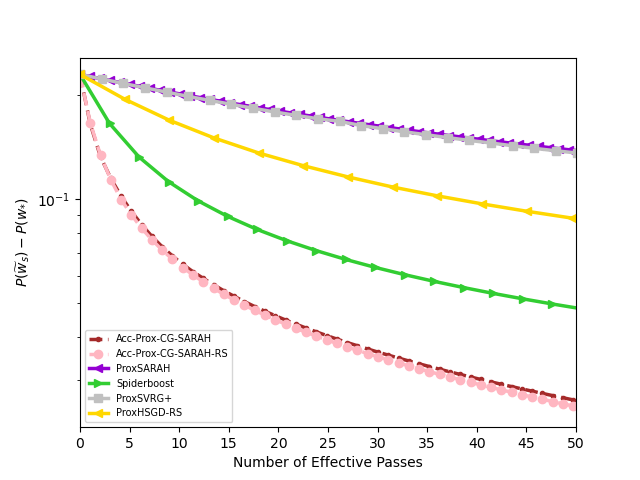}}
	
	\subfigure[w8a]
	{\includegraphics[width=0.327\textwidth]{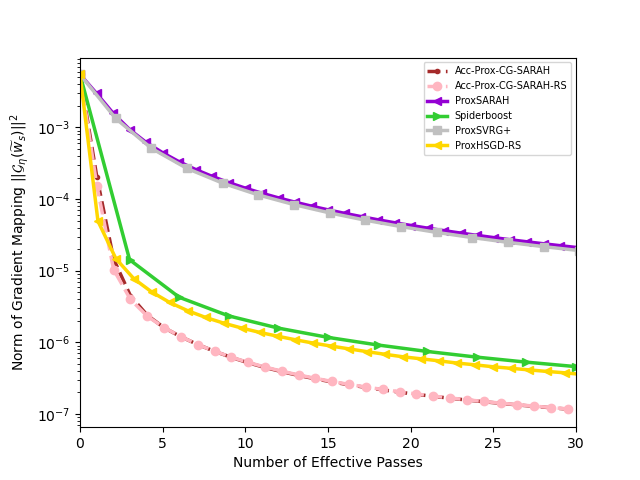}}
	\subfigure[a9a]
	{\includegraphics[width=0.327\textwidth]{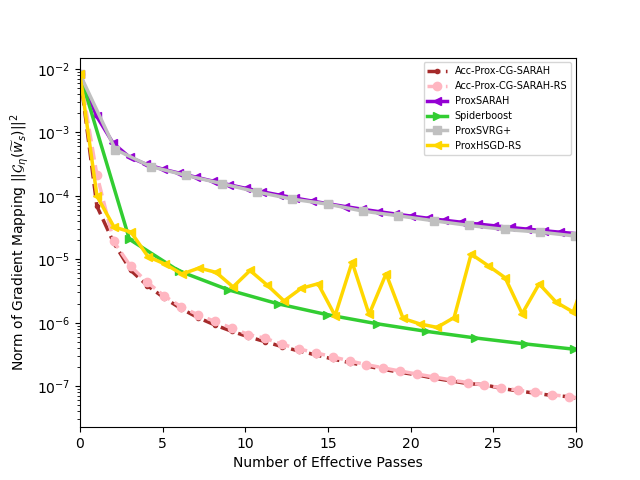}}
	\subfigure[gisette]
	{\includegraphics[width=0.327\textwidth]{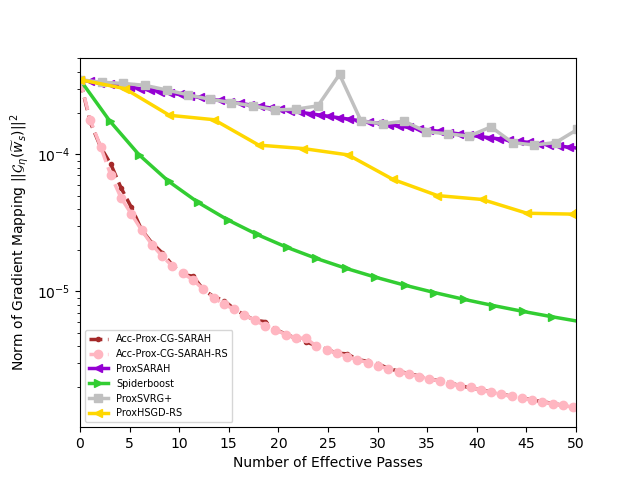}}
	\caption{\footnotesize Comparisons with other state-of-art methods over nonconvex loss $\ell_d$.}
	\label{fig7}
\end{figure*}

Fig. \ref{fig4} - Fig. \ref{fig7} illustrate the performance of these algorithms on all the three data sets. From Fig. \ref{fig4} we see that Algorithm \ref{alg1} and Algorithm \ref{alg2} outperform the other methods consistently in terms of lorenz loss $\ell_a$. Particularly in $a9a$, the optimization progress of the other four methods appear to stagnate after $2-3$ effective passes, while our two algorithms advance forward smoothly and reach notably small values.

Fig. \ref{fig5} presents the results over normalized sigmoid loss model with $\ell_b$. Algorithm \ref{alg1} and Algorithm \ref{alg2} consistently converge faster throughout the entire optimization in $w8a$ and $gisette$. As in $a9a$, ProxSARAH, ProxHSGD-RS and Spiderboost exhibit impressive performance in early periods, but are soon overtaken after approximately $3$ effective passes. Notice that Spiderboost yields the smallest gradient mapping norm.

Fig. \ref{fig6} reports the sub-optimality and the gradient mapping norms for the logistic difference loss model in terms of $\ell_c$. Our two algorithms are the most effective for every case and both converge to better levels in all the data sets. ProxSARAH remains inactive during the tests, which demonstrates the efficacy of conjugacy. In addtion, the gradient mapping norms of ProxHSGD-RS oscillate in $a9a$ and $gisette$.

Fig. \ref{fig7} shows the results associated with nonconvex loss $\ell_d$. Both Algorithm \ref{alg1} and Algorithm \ref{alg2} are clearly advantageous than the other competitors. Notably, the gradient mapping norms of ProxHSGD-RS oscillate severely in $w8a$.

Overall, Acc-Prox-CG-SARAH and Acc-Prox-CG-SARAH-RS have demonstrated significant competitiveness with modern state-of-art methods. ProxSARAH, Spiderboost, ProxSVRG+ and ProxHSGD-RS may be comparable in certain cases, but soon expose relative powerlessness on the others. Our frameworks can apply large adaptive step sizes due to the line search, which ensures fast and robust convergence process. Furthermore, we observe that the slight improvement from the deterministic restart scheme may be mitigated by the inherent stochasticity of optimization.

\subsection{Performance of the practical variant}
Then, we proceed with the experiments for our practical variant: Algorithm \ref{alg4}. We choose the switching frequency $t=5$ and settle the switching fixed step size $\eta_f = \frac{1}{L}$. Experimental results are displayed in Table \ref{table4} and Fig. \ref{fig7} (We select partial results to avoid overloading the paper).
\begin{table}[h]
	\setlength{\abovecaptionskip}{0pt}
	\setlength{\belowcaptionskip}{1pt}
	\centering
	\vspace{4pt}
	\renewcommand\arraystretch{1.5}
	\resizebox{0.5\textwidth}{!}{
		\begin{tabular}{cccc}
			\hline
			\multicolumn{1}{l}{Models/Data sets}&\multicolumn{1}{c}{$w8a$} &\multicolumn{1}{c}{$a9a$} &\multicolumn{1}{c}{$gisette$} 
			\\ \hline\noalign{\smallskip}
			$\ell_a$ & \makecell{\textcolor{orange}{37.6513}\\ \textcolor{red}{37.3152}\\ \textcolor{blue}{26.8021} }  & \makecell{\textcolor{orange}{5.7829}\\ \textcolor{red}{5.7244}\\ \textcolor{blue}{5.0187} } & \makecell{ \textcolor{orange}{35.9311}\\ \textcolor{red}{34.4444}\\ \textcolor{blue}{25.1567}} \\
			\hline
			$\ell_b$ &  \makecell{\textcolor{orange}{4.8118}\\ \textcolor{red}{4.6603}\\ \textcolor{blue}{5.2827} }  &  \makecell{\textcolor{orange}{7.0796}\\ \textcolor{red}{6.8312}\\ \textcolor{blue}{6.3045} } & \makecell{\textcolor{orange}{36.5994}\\ \textcolor{red}{36.5718}\\ \textcolor{blue}{33.6973}} \\
			\hline
			$\ell_c$  & \makecell{\textcolor{orange}{86.8253}\\ \textcolor{red}{84.4777}\\ \textcolor{blue}{79.9164} }  & \makecell{\textcolor{orange}{54.0298}\\ \textcolor{red}{53.8130}\\ \textcolor{blue}{46.5746} } & \makecell{\textcolor{orange}{65.0400}\\ \textcolor{red}{61.3776}\\ \textcolor{blue}{34.2004}} \\
			\hline
			$\ell_d$ & \makecell{\textcolor{orange}{66.0075}\\ \textcolor{red}{64.5754}\\ \textcolor{blue}{57.5822} }  & \makecell{\textcolor{orange}{270.7704}\\ \textcolor{red}{247.4300}\\ \textcolor{blue}{36.2046} } & \makecell{\textcolor{orange}{39.0883}\\ \textcolor{red}{39.0911}\\ \textcolor{blue}{34.4368} }  \\
			\hline
	\end{tabular}}
	\caption{The average running time (20 runs) of Algorithm \ref{alg1} (marked in \textcolor{orange}{orange}), Algorithm \ref{alg2} (marked in \textcolor{red}{red}) and Algorithm \ref{alg4} (marked in \textcolor{blue}{blue}) on different models. All were implemented in Python 3.9.5 version, with an AMD Ryzen 7 4800H processor 2.90 GHz and 16 GB of RAM.} 
	\label{table4}
\end{table}

\begin{figure*}[h]
	\centering
	\subfigure[w8a, $\ell_a$ ]
	{\includegraphics[width=0.327\textwidth]{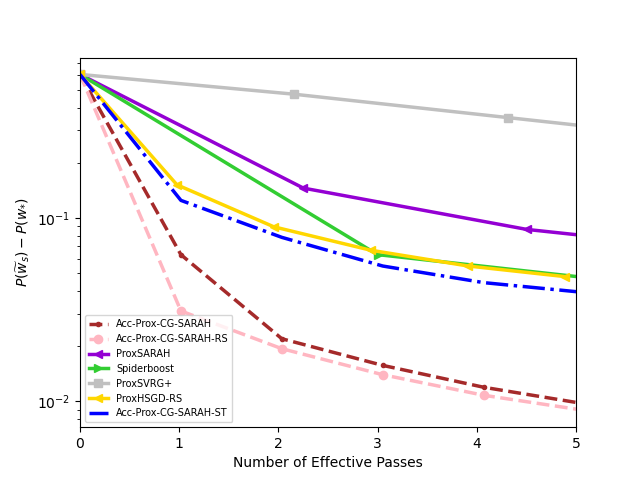}}
	\subfigure[w8a, $\ell_b$]
	{\includegraphics[width=0.327\textwidth]{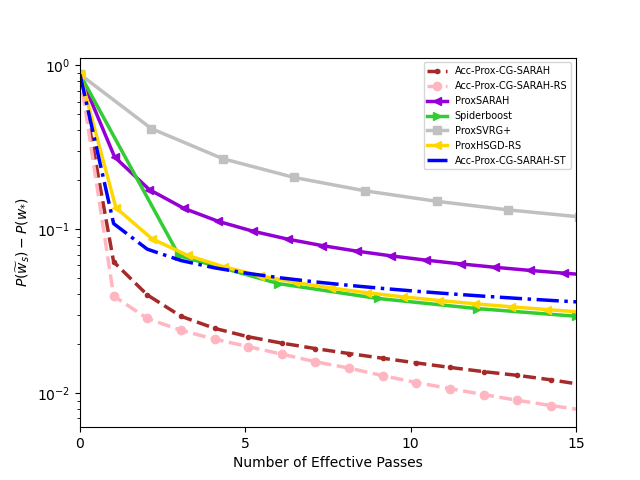}}
	\subfigure[a9a, $\ell_b$]
	{\includegraphics[width=0.327\textwidth]{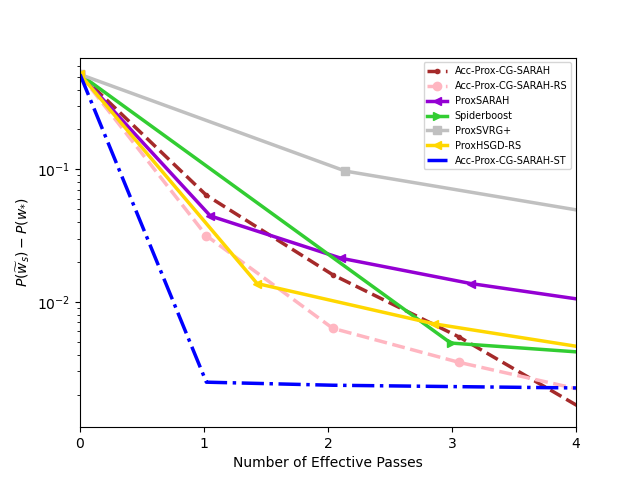}}
	
	\subfigure[w8a, $\ell_a$]
	{\includegraphics[width=0.327\textwidth]{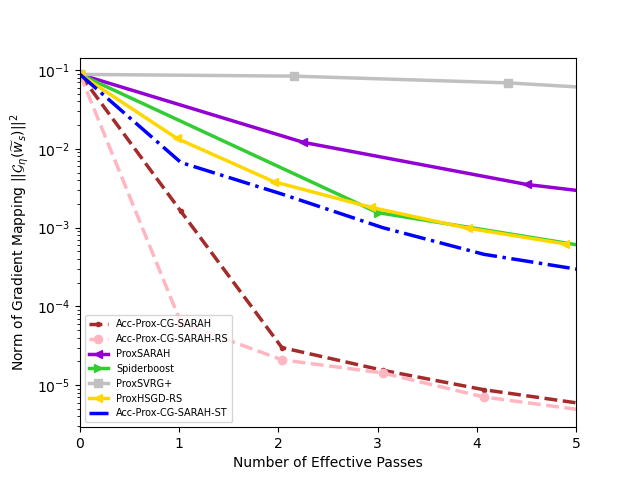}}
	\subfigure[w8a, $\ell_b$]
	{\includegraphics[width=0.327\textwidth]{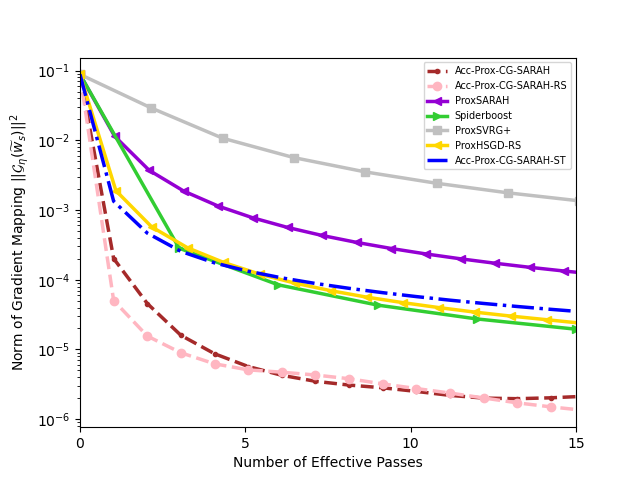}}
	\subfigure[a9a, $\ell_b$]
	{\includegraphics[width=0.327\textwidth]{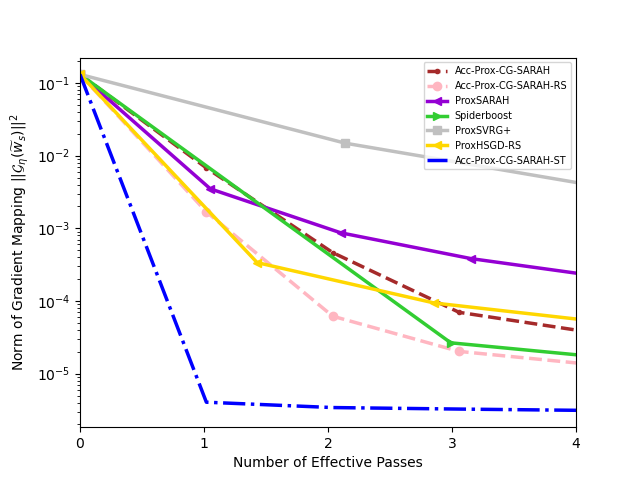}}
	\caption{\footnotesize Several representative results regarding our practical variant Algorithm \ref{alg4}.}
	\label{fig8}
\end{figure*}

From Table \ref{table4}, we observe that the average running time decreases as we switch to reduce the number of calls to line search. The average running time of Acc-Prox-CG-SARAH-RS is less than Acc-Prox-CG-SARAH-RS in general. Furthermore, as depicted in Fig. \ref{fig7}, our practical variant Acc-Prox-CG-SARAH-ST still matches the speed of modern state-of-art algorithms.

\section{Conclusions}
In the paper, we propose two new stochastic conjugate frameworks for a class of nonconvex and possibly nonsmooth optimization problems. By referring to the conventional conjugate gradient (CG) method, we devise a deterministic restart scheme for stochastic optimization and apply it in our second approach. Considering the computational cost, we further develop another practical variant.

Then, we rigorously prove that Acc-Prox-CG-SARAH, Acc-Prox-CG-SARAH-RS and Acc-Prox-CG-SARAH-ST all obtain linear convergence rates in nonconvex and nonsmooth optimization. In addition, we extend our analysis under the well-known gradient dominant condition.

Experiments have demonstrated that Acc-Prox-CG-SARAH and Acc-Prox-CG-SARAH-RS outperform state-of-art methods consistently and Acc-Prox-CG-SARAH-ST can achieves comparable convergence speed.

\bibliography{sn-bibliography}% common bib file

\newpage
%\appendix

\section*{\centerline{Appendix}}

\paragraph{A. Proof of Lemma \ref{lemma2}}
On the basis of Theorem 2 in \cite{jin2018cg-svrg}, we have 
\begin{equation}
\label{lemma2_1}
\langle-v_{k}^{(s)},d_{k-1}^{(s)}\rangle\leq-c_2\langle v_{k-1}^{(s)},d_{k-1}^{(s)}\rangle\leq\frac{c_2}{1-c_2}\|v_{k-1}^{(s)}\|^2.
\end{equation}
By Assumption 2, we obtain
\begin{equation}
\label{lemma2_2}
\mathbb{E}\left[\|v_{k}^{(s)}\|^2\right]\leq\hat{\beta} \mathbb{E}\left[\|v_{k-1}^{(s)}\|^2\right].
\end{equation}
Combining (\ref{lemma2_1}) and (\ref{lemma2_2}), it further implies
$$
\begin{aligned}
\mathbb{E}\left[\|d_k^{(s)}\|^{2}\right]&={\mathbb{E}\left[\|\beta_{k}d_{k-1}^{(s)}-v_{k}^{(s)}\|^{2}\right]}  \\
&=\mathbb{E}\left[\|v_{k}^{(s)}\|^{2}\right]-2\beta_{k}\mathbb{E}\left[\langle v_{k}^{(s)}, d_{k-1}^{(s)}\rangle\right]+\beta_k^{2}\mathbb{E}\left[\|d_{k-1}^{(s)}\|^{2}\right] \\
&\leq\hat{\beta}\mathbb{E}\left[\|v_{k-1}^{(s)}\|^2\right]+\frac{2\hat{\beta}c_2}{1-c_2}\mathbb{E}\left[\|v_{k-1}^{(s)}\|^2\right]+\hat{\beta}^{2}\mathbb{E}\left[\|d_{k-1}^{(s)}\|^{2}\right] \\
&=\hat{\beta}\frac{1+c_{2}}{1-c_{2}}\mathbb{E}\left[\|v_{k-1}^{(s)}\|^{2}\right]+\hat{\beta}^{2}\mathbb{E}\left[\|d_{k-1}^{(s)}\|^{2}\right].
\end{aligned}
$$
Since $\Phi(c_2)\leq\alpha$ ($\alpha>1$), we have
\begin{equation}
\label{lemma2_3}
\mathbb{E}\left[\|d_k^{(s)}\|^2\right]\leq\alpha\hat{\beta}\mathbb{E}\left[\|v_{k-1}^{(s)}\|^2\right]+\hat{\beta}^2\mathbb{E}\left[\|d_{k-1}^{(s)}\|^2\right].
\end{equation}
By applying (\ref{lemma2_2}) (\ref{lemma2_3}) recursively, we obtain
\begin{equation}
\label{Ed}
\begin{aligned}
&\mathbb{E}\left[\|d_{k}^{(s)}\|^{2}\right]\\
&\leq\alpha\hat{\beta}\left(\mathbb{E}\left[\|v_{k-1}^{(s)}\|^{2}\right]+\hat{\beta}^{2}\mathbb{E}\left[\|v_{k-2}^{(s)}\|^{2}\right]+\cdots+(\hat{\beta}^{2})^{k-1}\mathbb{E}\left[\|v_{0}^{(s)}\|^{2}\right]\right)+(\hat{\beta}^{2})^{k}\mathbb{E}\left[\|d_{0}^{(s)}\|^{2}\right] \\
&=\alpha\hat{\beta}\left(\hat{\beta}^{k-1}\mathbb{E}\left[\|v_{0}^{(s)}\|^{2}\right]+(\hat{\beta}^{2})\hat{\beta}^{k-2}\mathbb{E}\left[\|v_{0}^{(s)}\|^{2}\right]+\cdots+(\hat{\beta}^{2})^{k-1}\mathbb{E}\left[\|v_{0}^{(s)}\|^{2}\right]\right)+(\hat{\beta}^{2})^{k}\mathbb{E}\left[\|d_{0}^{(s)}\|^{2}\right] \\
&=\alpha\hat{\beta}^k\left(\sum_{j=0}^{k-1}\hat{\beta}^j\right)\mathbb{E}\left[\|v_0^{(s)}\|^2\right]+\hat{\beta}^{2k}\mathbb{E}\left[\|v_0^{(s)}\|^2\right]+\hat{\beta}^{2k}\mathbb{E}\left[\|d_{0}^{(s)}\|^{2}\right]-\hat{\beta}^{2k}\mathbb{E}\left[\|v_0^{(s)}\|^2\right] \\
&=\left(\alpha\hat{\beta}^{k}\frac{1-\hat{\beta}^{k}}{1-\hat{\beta}}+\hat{\beta}^{2k}\right)\mathbb{E}\left[\|v_0^{(s)}\|^2\right]+\hat{\beta}^{2k}\mathbb{E}\left[\|d_{0}^{(s)}\|^{2}-\|v_0^{(s)}\|^2\right] \\
&=\left(\frac{\alpha}{1-\hat{\beta}}\hat{\beta}^{k}-\frac{\alpha -1+\hat{\beta}}{1-\hat{\beta}}\hat{\beta}^{2k}\right)\mathbb{E}\left[\|v_0^{(s)}\|^2\right]+\hat{\beta}^{2k}\mathbb{E}\left[\left|\|d_{0}^{(s)}\|^{2}-\|v_0^s\|^{2}\right|\right]    
\end{aligned}
\end{equation}
Using (\ref{deviation}) in Assumption 3, we define $\mathbb{E}\left[\left|\|v_{k}^{(s)}\|^{2}-\|\nabla f(w_k^{(s)})\|^{2}\right|\right]\leq\hat{\sigma}^2$ for any $s, k\geq0$. Thus, there exists a constant $\tau>0$ such that $\hat{\sigma}^2=\tau\sigma^2$. In Algorithm \ref{alg1}, we have $d_{0}^{(s)}=v_m^{(s-1)}$ and $v_0^{(s)}=\nabla f(w_0^{(s)})=\nabla f(w_m^{(s-1)})$, then  
$\mathbb{E}\left[\|d_k^{(s)}\|^2\right]\leq g(k)\mathbb{E}\left[\|v_0^{(s)}\|^2\right]+\hat{\beta}^{2k}\tau\sigma^2$. In Algorithm \ref{alg2} with $d_{0}^{(s)}=-v_0^{(s)}=-\nabla f(w_0^{(s)})$, it's adequate to have $\mathbb{E}\left[\|d_k^{(s)}\|^2\right]\leq g(k)\mathbb{E}\left[\|v_0^{(s)}\|^2\right]$.

\vskip 0.5cm

\paragraph{B.  Proof of Lemma \ref{lemma3} }  
Since $w_{k+1} =(1-\gamma)w_{k}+\gamma y_{k}$, we have $w_{k+1}-w_{k}=\gamma(y_{k}-w_{k}).$
Due to the $L$-smoothness of $f$, we obtain by (\ref{L2}) that
\begin{equation}
\label{Lemma3_1}
\begin{aligned}
f(w_{k+1})& \leq f(w_{k})+\langle\nabla f(w_{k}),w_{k+1}-w_{k}\rangle+\frac{L}{2}\|w_{k+1}-w_{k}\|^{2}  \\
&=f(w_{k})+\gamma\langle\nabla f(w_{k}),y_{k}-w_{k}\rangle+\frac{L\gamma^{2}}{2}\|y_{k}-w_{k}\|^{2}.
\end{aligned}
\end{equation}
By   convexity of $\varphi$, we get
\begin{equation}
\label{Lemma3_2}
\varphi(w_{k+1})\le\left(1-\gamma\right)\varphi(w_{k})+\gamma\varphi(y_k)\le\varphi(w_{k})+\gamma \cdot \langle \widetilde{\partial}\varphi(y_k),y_k-w_{k}\rangle,
\end{equation}
where $\widetilde{\partial}\varphi$ denotes the subgradients of $\varphi$ and each $\widetilde{\partial}\varphi(y_k)\in\partial\varphi(y_{k})$.
Next, the optimality condition of $\mathrm{prox}_{\eta\varphi}(w_{k}+\eta d_{k})$ implies 
\begin{equation}
\label{Lemma3_4}
\widetilde{\partial}\psi(y_{k})=d_{k}-\frac{1}{\eta}(y_{k}-w_{k})
\end{equation}
for some $\widetilde{\partial}\psi(y_{k})\in\partial\psi(y_{k})$. Substituting (\ref{Lemma3_4}) into (\ref{Lemma3_2}), we attain
\begin{equation}
\label{Lemma3_5}
\varphi(w_{k+1})\leq\varphi(w_{k})+\gamma\langle d_{k},y_{k}-w_{k}\rangle-\frac{\gamma}{\eta}\|y_{k}-w_{k}\|^{2}.
\end{equation}
Combining (\ref{Lemma3_1}) (\ref{Lemma3_5}), and   using $\langle a,b\rangle={\frac{1}{2}}[\|a\|^{2}+\|b\|^{2}-\|a-b\|^{2}]$ yields
\begin{equation}
\label{Lemma3_6}
\begin{aligned}
&P(w_{k+1})\\
&\leq P(w_{k})+\gamma\langle\nabla f(w_{k})+d_{k},y_{k}-w_{k}\rangle-\left(\frac{\gamma}{\eta}-\frac{L\gamma^{2}}{2}\right)\|y_{k}-w_{k}\|^{2}\\
&=P(w_{k})+\frac{\gamma}{2}\|\nabla f(w_{k})+d_{k}\|^{2}+\left(\frac{\gamma}{2}-\frac{\gamma}{\eta}+\frac{L\gamma^{2}}{2}\right)\|y_{k}-w_{k}\|^{2}-\Lambda_{k},
\end{aligned}
\end{equation}
where $\Lambda_{k}=\frac{\gamma}{2}\|\nabla f(w_{k})+d_{k}-(y_{k}-w_{k})\|^{2}\geq0$.

Note from (\ref{mapping}) that $\mathcal{G}_\eta(w)=\frac{1}{\eta}\left(w-\operatorname{prox}_{\eta \varphi}\left(w-\eta \nabla f(w)\right)\right)$ is the gradient mapping of $P$, thus  
\begin{equation}
\label{Lemma3_7}
\eta\|\mathcal{G}_{\eta}(w_{k})\|=\|w_{k}-\mathrm{prox}_{\eta\varphi}(w_{k}-\eta\nabla f(w_{k}))\|.
\end{equation}
Using the contractiveness of $\|\operatorname{prox}_{\eta\varphi}(x)-\operatorname{prox}_{\eta \varphi}(z)\|\leq\|x-z\|$ (Moreau and Jean Jacques, 1962 \cite{moreau1962fonctions}) and the triangle inequality, we derive from (\ref{Lemma3_7}) that
$$
\begin{aligned}
&\eta\|\mathcal{G}_{\eta}(w_{k})\|\\
&\leq\|y_{k}-w_{k}\|+\|\mathrm{prox}_{\eta\varphi}(w_{k}-\eta\nabla f(w_{k}))-y_{k}\|  \\
&=\|y_{k}-w_{k}\|+\|\mathrm{prox}_{\eta\varphi}(w_{k}-\eta\nabla f(w_{k}))-\mathrm{prox}_{\eta\varphi}(w_{k}+\eta d_{k})\| \\
&\leq\|y_{k}-w_{k}\|+\eta\|\nabla f(w_{k})+d_{k}\|.
\end{aligned}
$$
Due to $\|a+b\|^2\leq2\|a\|^2+2\|b\|^2$, we further have
\begin{equation}
\label{Lemma3_8}
\eta^2\|\mathcal{G}_{\eta}(w_k)\|^2\leq2\|y_{k}-w_k\|^2+2\eta^2\|\nabla f(w_k)+d_k\|^2.
\end{equation}

Now, multiplying (\ref{Lemma3_8}) with $\frac{\gamma}{2}>0$ and adding it to (\ref{Lemma3_6}) provides
$$
\begin{aligned}
P(w_{k+1}) &\leq P(w_{k})-\frac{\gamma\eta^{2}}{2}\|\mathcal{G}_{\eta}(w_{k})\|^{2}+\frac{\gamma}{2}\left(1+2\eta^{2}\right)\|\nabla f(w_{k})+d_{k}\|^{2} \\
&-\frac{\gamma}{2}\left(\frac{2}{\eta}-L\gamma-3\right)\|y_{k}-w_{k}\|^{2}-\Lambda_{k}.
\end{aligned}
$$
By considering $\eta\in\{\eta_1, \eta_2\}$, we take the full expectation to obtain
\begin{equation}
\label{Lemma3_9}
\begin{aligned}
\mathbb{E}\left[ P(w_{k+1})\right] &\leq \mathbb{E}\left[P(w_{k})\right]-\frac{\gamma\eta_1^{2}}{2}\mathbb{E}\left[\|\mathcal{G}_{\eta}(w_{k})\|^{2}\right]+\frac{\gamma}{2}\left(1+2\eta_2^{2}\right)\mathbb{E}\left[\|\nabla f(w_{k})+d_{k}\|^{2}\right] \\
&-\frac{\gamma}{2}\left(\frac{2}{\eta_2}-L\gamma-3\right)\mathbb{E}\left[\|y_{k}-w_{k}\|^{2}\right]-\mathbb{E}\left[\Lambda_{k}\right].
\end{aligned}
\end{equation}
For any $k\geq1$ in (\ref{Lemma3_9}), it holds that
$$
\begin{aligned}
&\mathbb{E}\left[P(w_{k+1})\right]\\
&\leq \mathbb{E}\left[P(w_{k})\right]-\frac{\gamma\eta_1^{2}}{2}\mathbb{E}\left[\|\mathcal{G}_{\eta}(w_{k})\|^{2}\right]+\gamma\left(1+2\eta_2^{2}\right)\mathbb{E}\left[\|\nabla f(w_{k})-v_{k}\|^{2}\right] \\
&+\gamma\left(1+2\eta_2^{2}\right)\hat{\beta}^2\mathbb{E}\left[\|d_{k-1}\|^2\right]-\frac{\gamma}{2}\left(\frac{2}{\eta_2}-L\gamma-3\right)\mathbb{E}\left[\|y_{k}-w_{k}\|^{2}\right]-\mathbb{E}\left[\Lambda_{k}\right].
\end{aligned}
$$
Summing over it by $k=1,\ldots,m$, we have
\begin{equation}
\label{Lemma3_10}
\begin{aligned}
&\mathbb{E}\left[P(w_{m+1})\right]\\
&\leq \mathbb{E}\left[P(w_{1})\right]+\gamma\left(1+2\eta_2^{2}\right)\sum_{k=1}^{m}\mathbb{E}\left[\|\nabla f(w_{k})-v_{k}\|^{2}\right]+ \gamma\left(1+2\eta_2^{2}\right)\hat{\beta}^2\sum_{k=1}^{m}\mathbb{E}\left[\|d_{k-1}\|\right]^2 \\
&-\frac{\gamma\eta_1^{2}}{2}\sum_{k=1}^{m}\mathbb{E}\left[\|\mathcal{G}_{\eta}(w_{k})\|^{2}\right]-\frac{\gamma}{2}\left(\frac2{\eta_2}-L\gamma-3\right)\sum_{k=1}^{m}\mathbb{E}\left[\|y_{k}-w_{k}\|^{2}\right]-\sum_{k=1}^{m}\mathbb{E}\left[\Lambda_{k}\right].  
\end{aligned}
\end{equation}
For $k=0$ in (\ref{Lemma3_9}), we point out that $d_0=-\nabla f(w_0)=-v_0$, which implies
\begin{equation}
\label{Lemma3_11}
\begin{aligned}
\mathbb{E}\left[P(w_1)\right]
&\leq \mathbb{E}\left[P(w_{0})\right]+\gamma\left(1+2\eta_2^{2}\right)\mathbb{E}\left[\|\nabla f(w_{0})-v_{0}\|^{2}\right] \\
&-\frac{\gamma\eta_1^{2}}{2}\mathbb{E}\left[\|\mathcal{G}_{\eta}(w_{0})\|^{2}\right]-\frac{\gamma}{2}\left(\frac2{\eta_2}-L\gamma-3\right)\mathbb{E}\left[\|y_{1}-w_{0}\|^{2}\right]-\mathbb{E}\left[\Lambda_{0}\right]. 
\end{aligned}
\end{equation}
After integrating (\ref{Lemma3_10}) (\ref{Lemma3_11}) and dropping $-\sum_{k=0}^{m}\mathbb{E}\left[\Lambda_{k}^{(s)}\right]\leq0$, we finally get
$$
\begin{aligned}
&\mathbb{E}\left[P(w_{m+1})\right]\\
&\leq \mathbb{E}\left[P(w_{0})\right]-\frac{\gamma\eta_1^{2}}{2}\sum_{k=0}^{m}\mathbb{E}\left[\|\mathcal{G}_{\eta}(w_{k})\|^{2}\right]+ \gamma\left(1+2\eta_2^{2}\right)\hat{\beta}^2\sum_{k=1}^{m}\mathbb{E}\left[\|d_{k-1}\|^2\right] \\
&+\gamma\left(1+2\eta_2^{2}\right)\sum_{k=0}^{m}\mathbb{E}\left[\|\nabla f(w_{k})-v_{k}\|^{2}\right]-\frac{\gamma}{2}\left(\frac2{\eta_2}-L\gamma-3\right)\sum_{k=0}^{m}\mathbb{E}\left[\|y_{k}-w_{k}\|^{2}\right].
\end{aligned}
$$

\vskip 0.5cm

\paragraph{C.  Proof of Theorem \ref{theorem1} }

%\section{Appendix: the proof of Theorem \ref{theorem1}}
Firstly, we proceed with the proof for Algorithm \ref{alg1}. Since $d_0^{(s)}=-v_m^{(s-1)}$ and $w_0^{(s)}=w_{m}^{(s-1)}$ for all $s>1$, by slightly modifying (\ref{Lemma3_10}) in Lemma \ref{lemma3} we have
\begin{equation}
\label{Thm_1}
\begin{aligned}
&\mathbb{E}\left[P(w_{m+1}^{(s)})\right]\\
&\leq \mathbb{E}\left[P(w_{0}^{(s)})\right]+\gamma\left(1+2\eta_2^{2}\right)\sum_{k=0}^{m}\mathbb{E}\left[\|\nabla f(w_{k}^{(s)})-v_{k}^{(s)}\|^{2}\right]+ \gamma\left(1+2\eta_2^{2}\right)\hat{\beta}^2\sum_{k=1}^{m}\mathbb{E}\left[\|d_{k-1}^{(s)}\|^2\right] \\
&-\frac{\gamma\eta_1^{2}}{2}\sum_{k=0}^{m}\mathbb{E}\left[\|\mathcal{G}_{\eta}(w_{k}^{(s)})\|^{2}\right]-\frac{\gamma}{2}(\frac2{\eta_2}-L\gamma-3)\sum_{k=0}^{m}\mathbb{E}\left[\|y^{(s)}_{k}-w_{k}^{(s)}\|^{2}\right]+\frac{\gamma}{2}\left(1+2\eta_2^2\right)\sigma^2.
\end{aligned}
\end{equation}
According to Lemma 2 in \cite{nguyen2017MBSARAH}, we obtain
$$\mathbb{E}\left[\|\nabla f(w_k^{(s)})-v_{k}^{(s)}\|^2\right]=\sum_{j=1}^k\mathbb{E}\left[\|v_j^{(s)}-v_{j-1}^{(s)}\|^2\right]-\sum_{j=1}^k\mathbb{E}\left[\|\nabla f(w_j^{(s)})-\nabla f(w_{j-1}^{(s)})\|^2\right].$$
Based on Lemma \ref{lemma1} and (\ref{L3}), for all $j\geq1$, we have
$$
\begin{aligned}
\mathbb{E}\left[\|v_{j}^{(s)}-v_{j-1}^{(s)}\|^{2}\right]-\mathbb{E}\left[\|\nabla f(w_{j}^{(s)})-\nabla f(w_{j-1}^{(s)})\|^{2}\right]\leq\frac{\left(n-b\right)L^{2}\gamma^{2}}{b\left(n-1\right)}\|y_{j-1}^{(s)}-w_{j-1}^{(s)}\|^{2}.
\end{aligned}
$$
Rearranging the results yields
\begin{equation}
\label{Thm_2}
\sum_{k=0}^{m}\mathbb{E}\left[\|\nabla f(w_{k}^{(s)})-v_{k}^{(s)}\|^{2}\right]\leq\frac{\left(n-b\right)L^{2}\gamma^{2}}{b\left(n-1\right)}\sum_{k=1}^{m}\sum_{j=1}^{k}\|y_{j-1}^{(s)}-w_{j-1}^{(s)}\|^{2}.
\end{equation}
Furthermore, by using Lemma \ref{lemma2} we conclude that
\begin{equation}
\begin{aligned}
\sum_{k=1}^{m}\mathbb{E}\left[\|d_{k-1}^{(s)}\|^2\right]&=\sum_{k=0}^{m-1}\mathbb{E}\left[\|d_{k}^{(s)}\|^2\right]\\
&\leq\mathbb{E}\left[\|v_0^{(s)}\|^2\right]\sum_{k=0}^{m-1}g(k)+\tau\sigma^2\sum_{k=0}^{m-1}\hat{\beta}^{2k}\\
&=\frac{\left((\alpha-1)\hat{\beta}+1-(\alpha-1+\hat{\beta})\hat{\beta}^m\right)(1-\hat{\beta}^m)}{(1-\hat{\beta})^2(1+\hat{\beta})}\mathbb{E}\left[\|\nabla f(w_0^{(s)})\|^2\right]+\frac{1-\hat{\beta}^{2m}}{1-\hat{\beta}^2}\tau\sigma^2\\
&\leq\frac{\alpha}{(1-\hat{\beta})^2(1+\hat{\beta})}\mathbb{E}\left[\|\nabla f(w_0^{(s)})\|^2\right]+\frac{1}{1-\hat{\beta}^2}\tau\sigma^2.\label{Thm_3}
\end{aligned}
\end{equation}
Plugging (\ref{Thm_2}) (\ref{Thm_3}) into (\ref{Thm_1}) and applying $v_0^{(s)}=\nabla f(w_0^{(s)})$, it holds that
\begin{equation}
\label{Thm_4}
\begin{aligned}
&\mathbb{E}\left[P(w_{m+1}^{(s)})\right]\\
&\leq \mathbb{E}\left[P(w_{0}^{(s)})\right]-\frac{\gamma\eta_1^{2}}{2}\sum_{k=0}^{m}\mathbb{E}\left[\|\mathcal{G}_{\eta}(w_{k}^{(s)})\|^{2}\right]-\frac{\gamma}{2}(\frac2{\eta_2}-L\gamma-3)\sum_{k=0}^{m}\mathbb{E}\left[\|y^{(s)}_{k}-w_{k}^{(s)}\|^{2}\right] \\
&+\left(1+2\eta_2^2\right)\frac{\left(n-b\right)L^{2}\gamma^{3}}{b\left(n-1\right)}\sum_{k=1}^{m}\sum_{j=1}^{k}\mathbb{E}\left[\|y_{j}^{(s)}-w_{j-1}^{(s)}\|^{2}\right]+\gamma\left(1+2\eta_2^{2}\right)(\frac{\hat{\beta}^2}{1-\hat{\beta}^2}\tau+\frac{1}{2})\sigma^2\\
&+\gamma\left(1+2\eta_2^{2}\right)\frac{\alpha\hat{\beta}^2}{(1-\hat{\beta})^2(1+\hat{\beta})}\mathbb{E}\left[\|\nabla f(w_0^{(s)})\|^2\right].
\end{aligned}
\end{equation}
Subsequently, one can show that
$$
\begin{aligned}
&\sum_{k=1}^{m}\sum_{j=1}^{k}\mathbb{E}\left[\|y_{j-1}^{(s)}-w_{j-1}^{(s)}\|^{2}\right]\\
&=m\cdot\mathbb{E}\left[\|y_{0}^{(s)}-w_{0}^{(s)}\|^{2}\right]+\left(m-1\right)\cdot\mathbb{E}\left[\|y_{1}^{(s)}-w_{1}^{(s)}\|^{2}\right]+\cdots+ \mathbb{E}\left[\|y_{m-1}^{(s)}-w_{m-1}^{(s)}\|^{2}\right]\\
&\leq m\sum_{k=0}^{m}\mathbb{E}\left[\|y^{(s)}_{k}-w_{k}^{(s)}\|^{2}\right]
\end{aligned}
$$

Choosing $b$, $\gamma$, which satisfies %such that
\begin{equation}
	\label{Thm1}
	\left(2+4\eta_2^2\right)\frac{n-b}{b(n-1)}L^2\gamma^2m-(\frac{2}{\eta_2}-L\gamma-3)\leq0.
\end{equation}
Then, by (\ref{Thm1}) we have 
$$
\begin{aligned}
&-\frac{\gamma}{2}(\frac2{\eta_2}-L\gamma-3)\sum_{k=0}^{m}\mathbb{E}\left[\|y^{(s)}_{k}-w_{k}^{(s)}\|^{2}\right]+\left(1+2\eta_2^2\right)\frac{\left(n-b\right)L^{2}\gamma^{3}}{b\left(n-1\right)}\sum_{k=1}^{m}\sum_{j=1}^{k}\mathbb{E}\left[\|y_{j-1}^{(s)}-w_{j-1}^{(s)}\|^{2}\right]\\
&\leq\left(-\frac{\gamma}{2}(\frac2{\eta_2}-L\gamma-3)+\left(1+2\eta_2^2\right)\frac{\left(n-b\right)L^{2}\gamma^{3}m}{b\left(n-1\right)}\right)\sum_{k=0}^{m}\mathbb{E}\left[\|y^{(s)}_{k}-w_{k}^{(s)}\|^{2}\right]\\
&\leq0.
\end{aligned}
$$
Thus 
\begin{equation}
\begin{aligned}
&\mathbb{E}\left[P(w_{m+1}^{(s)})\right]\\
\leq  & \mathbb{E}\left[P(w_{0}^{(s)})\right]-\frac{\gamma\eta_1^{2}}{2}\sum_{k=0}^{m}\mathbb{E}\left[\|\mathcal{G}_{\eta}(w_{k}^{(s)})\|^{2}\right]\\
& +\gamma\left(1+2\eta_2^{2}\right)\frac{\alpha\hat{\beta}^2}{(1-\hat{\beta})^2(1+\hat{\beta})}\mathbb{E}\left[\|\nabla f(w_0^{(s)})\|^2\right]\\
&+\gamma\left(1+2\eta_2^{2}\right)\left(\frac{\hat{\beta}^2}{1-\hat{\beta}^2}\tau+\frac{1}{2}\right)\sigma^2.
\end{aligned}
\end{equation}
Since $w_{\star}=\arg\operatorname*{min}_{w}P(w)$ and ${\widetilde{w}}_{s}\sim U(\{w_{k}^{(s)}\})$, we obtain
$$
\begin{aligned}
&\mathbb{E}\left[\|\mathcal{G}_{\eta}(w_{m}^{(s)})\|^{2}\right]\\
\leq &\frac{1}{m+1}\left\{\frac{2}{\gamma\eta_1^2}\mathbb{E}\left[P(w_{0}^{(s)})-P(w_{\star})\right]\right. \\ 
& +\left.\frac{2+4\eta_2^2}{\eta_1^2}\left(\frac{\alpha\hat{\beta}^2}{(1-\hat{\beta})^2(1+\hat{\beta})}\mathbb{E}\left[\|\nabla f(w_0^{(s)})\|^2\right]+(\frac{\hat{\beta}^2}{1-\hat{\beta}^2}\tau+\frac{1}{2})\sigma^2\right)\right\} \\
= &\frac{\xi(1-\hat{\beta})^2(1+\hat{\beta})}{\gamma\alpha(1+2\eta_2^2)\hat{\beta}^2}\mathbb{E}\left[P(w_{0}^{(s)})-P(w_{\star})\right]+\xi\mathbb{E}\left[\|\nabla f(w_0^{(s)})\|^2\right]+C\xi,
\end{aligned}
$$
where we employ the definitions of $\xi$ and $C$ (\ref{C}) in the second equality.

\vspace{4pt}
For Algorithm \ref{alg2}, we employ the second result in Lemma \ref{lemma2} and smoothly obtain:
$$
\mathbb{E}\left[\|\mathcal{G}_{\eta}(w_{m}^{(s)})\|^{2}\right]
\leq\frac{\xi(1-\hat{\beta})^2(1+\hat{\beta})}{\gamma\alpha(1+2\eta_2^2)\hat{\beta}^2}\mathbb{E}\left[P(w_{0}^{(s)})-P(w_{\star})\right]+\xi\mathbb{E}\left[\|\nabla f(w_0^{(s)})\|^2\right].
$$

\vskip 0.5cm

\paragraph{D.  Proof of Theorem \ref{theorem3} }
%\section{Appendix: the proof of Theorem \ref{theorem3}}

Based on the proof in Theorem \ref{theorem1}, it's adequate to obtain
$$
\begin{aligned}
&\mathbb{E}\left[P(w_{m+1}^{(s)})\right]\\
&\leq \mathbb{E}\left[P(w_{0}^{(s)})\right]+\gamma\left(1+2\eta_2^{2}\right)\sum_{k\in H}\mathbb{E}\left[\|\nabla f(w_{k}^{(s)})-v_{k}^{(s)}\|^{2}\right]+ \gamma\left(1+2\eta_2^{2}\right)\hat{\beta}^2\sum_{k\in H}\mathbb{E}\left[\|d_{k-t}^{(s)}\|^2\right] \\
&+\frac{\gamma\left(1+2\eta_2^{2}\right)}{2}\sum_{k\in H'/\{0\}}\mathbb{E}\left[\|\nabla f(w_{k}^{(s)})-v_{k}^{(s)}\|^{2}\right]-\frac{\gamma\eta_1^{2}}{2}\sum_{k=0}^{m}\mathbb{E}\left[\|\mathcal{G}_{\eta}(w_{k}^{(s)})\|^{2}\right]\\
&-\frac{\gamma}{2}(\frac2{\eta_2}-L\gamma-3)\sum_{k=0}^{m}\mathbb{E}\left[\|y^{(s)}_{k}-w_{k}^{(s)}\|^{2}\right]+\frac{\gamma}{2}\left(1+2\eta_2^2\right)\sigma^2.
\end{aligned}
$$
Following a similar line of reasoning, one can show that
$$
\begin{aligned}
&\mathbb{E}\left[P(w_{m+1}^{(s)})\right]\\
&\leq \mathbb{E}\left[P(w_{0}^{(s)})\right]+\left(1+2\eta_2^2\right)\frac{q\left(n-b\right)L^{2}\gamma^{3}}{b\left(n-1\right)}\sum_{k=1}^{m}\mathbb{E}\left[\|y_{k-1}^{(s)}-w_{k-1}^{(s)}\|^{2}\right]+\gamma\left(1+2\eta_2^{2}\right)\hat{\beta}^2\sum_{k\in H}\mathbb{E}\left[\|d_{k-t}^{(s)}\|^2\right] \\
&+\left(1+2\eta_2^2\right)\frac{(m-1-q)\left(n-b\right)L^{2}\gamma^{3}}{2b\left(n-1\right)}\sum_{k=1}^{m}\mathbb{E}\left[\|y_{k-1}^{(s)}-w_{k-1}^{(s)}\|^{2}\right]-\frac{\gamma\eta_1^{2}}{2}\sum_{k=0}^{m}\mathbb{E}\left[\|\mathcal{G}_{\eta}(w_{k}^{(s)})\|^{2}\right]\\
&-\frac{\gamma}{2}(\frac2{\eta_2}-L\gamma-3)\sum_{k=0}^{m}\mathbb{E}\left[\|y^{(s)}_{k}-w_{k}^{(s)}\|^{2}\right]+\frac{\gamma}{2}\left(1+2\eta_2^2\right)\sigma^2.
\end{aligned}
$$

Jorge et al., 2006 \cite{jorge2006numerical} provides an original proof for Lemma 5.6 (Page 125) in Fletcher-Reeves method by mathematical induction. Since we utilize $d_k^{(s)}=-v_{k}^{(s)}+\beta_kd_{k-t}^{(s)}$, employ the updates (\ref{AFR2}) or (\ref{FRPR2}), and perform the line search to ensure the condition (\ref{Wolfe3}), we can apply the same mathematical induction at regular intervals here to derive the results. By analogous deduction, we obtain
$$
\begin{aligned}
&\mathbb{E}\left[P(w_{m+1}^{(s)})\right]\\
&\leq \mathbb{E}\left[P(w_{0}^{(s)})\right]-\frac{\gamma\eta_1^{2}}{2}\sum_{k=0}^{m}\mathbb{E}\left[\|\mathcal{G}_{\eta}(w_{k}^{(s)})\|^{2}\right]-\frac{\gamma}{2}(\frac2{\eta_2}-L\gamma-3)\sum_{k=0}^{m}\mathbb{E}\left[\|y^{(s)}_{k}-w_{k}^{(s)}\|^{2}\right] \\
&+\left(1+2\eta_2^2\right)\frac{\left(n-b\right)L^{2}\gamma^{3}(m-1+q)}{2b\left(n-1\right)}\sum_{k=1}^{m}\mathbb{E}\left[\|y_{k-1}^{(s)}-w_{k-1}^{(s)}\|^{2}\right]+\gamma\left(1+2\eta_2^{2}\right)(\frac{\hat{\beta}^2(1-\hat{\beta}^{2q})}{1-\hat{\beta}^2}\tau+\frac{1}{2})\sigma^2\\
&+\gamma\left(1+2\eta_2^{2}\right)\frac{\alpha\hat{\beta}^2(1-\hat{\beta}^{q})}{(1-\hat{\beta})^2(1+\hat{\beta})}\mathbb{E}\left[\|\nabla f(w_0^{(s)})\|^2\right].
\end{aligned}
$$
Assume that we choose $b$, $\gamma$ such that ($1\leq q< m-1$)
	\begin{equation}
	\label{Thm3}
	\left(2+4\eta_2^2\right)\frac{n-b}{b(n-1)}L^2\gamma^2(m+q-1)-(\frac{2}{\eta_2}-L\gamma-3)\leq0.
	\end{equation}
By the definition of $\xi^{st}$ and $C^{st}$ (\ref{C2}) and utilizing condition (\ref{Thm3}), we obtain
$$
\mathbb{E}\left[\|\mathcal{G}_{\eta}(w_{m}^{(s)})\|^{2}\right]
\leq\frac{\xi^{st}(1-\hat{\beta})^2(1+\hat{\beta})}{\gamma\alpha(1+2\eta_2^2)\hat{\beta}^2(1-\hat{\beta}^{q})}\mathbb{E}\left[P(w_{0}^{(s)})-P(w_{\star})\right]+\xi^{st}\mathbb{E}\left[\|\nabla f(w_0^{(s)})\|^2\right]+C^{st}\xi^{st},
$$
where
$$
\xi^{st}=\frac{(2+4\eta_2^2)\alpha\hat{\beta}^2(1-\hat{\beta}^{q})}{(m+1)\eta_1^2(1-\hat{\beta})^2(1+\hat{\beta})},\quad
C^{st}=\left(\frac{(1-\hat{\beta})(1+\hat{\beta}^{q})}{\alpha}\tau+\frac{(1-\hat{\beta})^2(1+\hat{\beta})}{2\alpha\hat{\beta}^2(1-\hat{\beta}^{q})}\right)\sigma^2.
$$

\end{document}